\documentclass[11pt, reqno]{amsart}
\usepackage[utf8]{inputenc}
\usepackage{amsmath, amssymb, amsfonts, latexsym, amsthm}
\usepackage{calligra}
\usepackage{tcolorbox}

\newtheorem{thm}{Theorem}[section]
\newtheorem{cor}[thm]{Corollary}
\newtheorem{lem}[thm]{Lemma}
\newtheorem{prop}[thm]{Proposition}
\newtheorem{rem}[thm]{Remark}

\newtheorem{defn}[thm]{Definition}
\newtheorem{exam}[thm]{Example}
\numberwithin{equation}{section}

\newcommand{\mbb}{\mathbb}
\newcommand{\ra}{\rightarrow}

\newcommand{\pa}{\partial}
\newcommand{\ov}{\overline}

\newcommand{\Om}{\Omega}

\newcommand{\om}{\omega}

\newcommand{\D}{\mathbb D}

\newcommand{\G}{\mathbb G_2}

\begin{document}

	\title[Isometric embeddings of the unit disk into the symmetrized bidisk]{Holomorphic isometric embeddings of the unit disk into the symmetrized bidisk}
	\author{Aakanksha Jain} 
	\address{AJ: Department of Mathematics, Indian Institute of Science Education and Research, Pune 411008, India}
	\email{mathcv.jain@gmail.com, aakanksha.jain@iiserpune.ac.in}

\begin{abstract}
We give a complete characterization of the holomorphic isometric embeddings $F:(\mathbb{D},k\,ds^2)\rightarrow (\mathbb{G}_2,ds^2)$ of the unit disk into the symmetrized bidisk with respect to the Bergman metric up to a scalar constant. We show that the isometric constant $k$ can take only the values $k=1$ and $k=3$. For $k=3$, we prove that $(2z,z^2)$ is the unique isometry up to a reparametrization. For $k=1$, we construct an isometric embedding which intersects the royal variety at a single point and does not admit a holomorphic lift to the bidisk via the symmetrization map. We further show that this isometric embedding is unique up to reparametrizations. Finally, none of the isometries are totally geodesic.
\end{abstract}
	
\maketitle
	
\section{Introduction}

A strong rigidity phenomenon occurs for holomorphic isometric embeddings of an irreducible bounded symmetric domain $\Om$ into its product spaces $\Om^n$, $n\geq 2$, with respect to the Bergman metric up to a scalar constant. Motivated by a problem in number theory, Clozel--Ullmo \cite{CU} were led to study these isometries. Using the results of Mok \cite{hm2, hm1} on Hermitian metric rigidity, they explained that when $\operatorname{rank}(\Om)\geq 2$, any such isometric embedding must be totally geodesic. When $\operatorname{rank}(\Om) = 1$, that is, $\Om = \mathbb{B}^d$, Mok \cite{Mok1} proved that the embedding must also be totally geodesic for $d\geq 2$. However, when $\Om$ is the unit disk $\mathbb{D}$, Mok \cite{Mok} constructed a non-totally geodesic holomorphic isometric embedding into the polydisk $\mathbb{D}^n$, for each $n\geq 2$. Isometric embeddings and their characterizations have since been studied under different paradigms; see \cite{Chan, CXY, JV, Ng}, and the references therein. In particular, the construction of non-standard (non-totally geodesic) holomorphic isometric embeddings from the unit disk into bounded domains is an interesting problem.

\medskip

This manuscript studies holomorphic isometric embeddings from the unit disk into the symmetrized bidisk with respect to the Bergman metric up to scalar constants. The \textit{symmetrized bidisk}
\[
 \G=\{(\lambda_1+\lambda_2,\lambda_1\lambda_2):\lambda_1,\lambda_2\in\D\}
\]
is a non-homogeneous domain in several complex variables. The study of this domain was initiated by J. Agler, F. B. Yeh, and N. J. Young \cite{Agler_2, Agler_1, Agler-Young}. They were led to study the domain while attempting to solve particular cases of the $\mu$-synthesis problem in control theory. The domain has since been studied from the viewpoint of geometric function theory and operator theory; see, among others, \cite{ALY, ALY_extremal, Agler_0, TB1, Costara, Po, Ed, EZ, JP, Pal, PZ}. 

\medskip

Given a bounded domain $\Om\subset \mathbb{C}^n$, the Bergman space $A^2(\Om)$ is the Hilbert space consisting of holomorphic functions on $\Om$ that are square-integrable with respect to the Lebesgue volume measure $dV$. The reproducing kernel of $A^2(\Om)$ is called the Bergman kernel $K_{\Om}(z,w)$ of $\Om$. We shall denote $K_{\Om}(z):=K_{\Om}(z,z)$ for simplicity. The Bergman kernel, which can be written as
\[
K_{\Om}(z) = \sup\{\vert f(z) \vert^2 : f \in A^2(\Om), \Vert f \Vert_{L^2(\Om)} = 1 \},
\]
is a positive plurisubharmonic function on $\Om$. The Bergman metric 
\[
ds^2_{\Om} = \sum_{i, j} \frac{\pa^2}{\pa z_i \pa \ov z_j} \log K_{\Om}(z) \; dz_i \otimes d \ov z_j
\]
is, therefore, well-defined on $\Om$. We shall drop the subscript $\Om$ when the domain is clear from the context. 

\medskip

Let $k$ be a positive integer, and
\[
F:=(S,P):(\mathbb D,k\,ds^{2})
\rightarrow
(\mathbb G_2,ds^{2}).
\]
be a holomorphic isometric embedding with respect to the \textit{isometric constant} $k$. 
Let $\pi:\mathbb D \times \mathbb{D} \to\mathbb G_2$ denote the \textit{symmetrization map} 
\[
\pi(\lambda_1,\lambda_2)
=(\lambda_1+\lambda_2,\lambda_1\lambda_2).
\]
The map $\pi$ is a proper holomorphic map of multiplicity two. Moreover, $\pi: \mathbb{D}^2\setminus\Delta \rightarrow\mathbb{G}_2\setminus \mathcal{R}$ is a holomorphic covering, where $\Delta = \{(\lambda,\lambda): \lambda\in\mathbb{D}\}$ is the diagonal of the bidisk and $\mathcal{R}=\{(2\lambda,\lambda^2): \lambda\in\mathbb{C}\}$ is the \textit{royal variety}.  

\begin{defn}
A holomorphic isometric embedding $F:=(S, P): (\mathbb{D}, k \,ds^2)\rightarrow (\mathbb{G}_2,ds^2)$ is called liftable if there exists a holomorphic map $G:=(g_1,g_2):\mathbb{D}\rightarrow \mathbb{D}\times \mathbb{D}$ such that
\[
F = \pi\circ G ,
\qquad
\text{that is, }
\qquad
S= g_1 + g_2
\quad \text{and} \quad
P = g_1 g_2.
\]
The isometric embedding is called non-liftable if there does not exist any such $G$.
\end{defn}

The aim of this article is to investigate the existence and structure of holomorphic isometric embeddings $F:(\mathbb D,k\,ds^{2}) \rightarrow
(\mathbb G_2,ds^{2})$. In particular, we determine the possible values of the isometric constant $k$, study the liftability of $F$ under the symmetrization map, and give characterization results. The central result of the article is as follows.

\begin{thm}[Main Theorem]
Let 
\[
F:(\mathbb D,k\,ds^{2})\rightarrow (\mathbb G_2,ds^{2})
\]
be a holomorphic isometric embedding, where the $k$ is a positive integer. Then, $k\in\{1,3\}$. Moreover, up to a reparametrization of $F$, the following classification holds:
\begin{enumerate}
\item If $k=1$, then $F$ is given by
\[
F_N(z)=
\left(
\sqrt{2}(\sqrt{z+1}-1),
1-\sqrt{z+1}
\right).
\]

\item If $k=3$, then $F$ is given by
\[
F_R(z)=(2z,z^2).
\]
\end{enumerate}
Neither $F_N$ nor $F_R$ is a totally geodesic isometric embedding. Furthermore, $F_N$ is a non-liftable isometry. 
\end{thm}

The proof of the main result proceeds through the following steps.

\begin{enumerate}
\item Our approach is based on a functional equation arising from the reproducing kernels associated with the Bergman metrics. We first show that $k\in\{1,2,3\}$ and relate these values to the boundary cluster set of the isometric embedding. More precisely, the value $k=1$ corresponds to a cluster set meeting the non-distinguished boundary of $\mathbb{G}_2$, while $k=3$ forces the cluster set to lie in the royal variety. The intermediate case $k=2$ would require the cluster set to be contained in the distinguished boundary away from the royal variety. 

\item We completely characterize holomorphic isometric embeddings for $k=3$. Namely, every such embedding is the royal isometry
\[
F_R(z)=(2z,z^2),\qquad z\in\mathbb D
\]
up to a reparametrization, see Theorem \ref{uni3}.

\item For $k=1$, we construct the non-royal holomorphic isometric embedding
\[
F_N(z)=
\left(
\sqrt{2}(\sqrt{z+1}-1),
1-\sqrt{z+1}
\right),
\qquad z\in\mathbb D.
\]
Its image intersects the royal variety at exactly one point $(0,0)$. Moreover, we prove that neither $F_R$ nor $F_N$ is a totally geodesic isometric embedding. 

\item We prove that there does not exist any liftable holomorphic isometric embedding intersecting the royal variety, for $k=1$ or $k=2$, see Theorem \ref{lift}. The key observation is that any isometric embedding intersecting the royal variety must do so transversely. This transversality obstructs the existence of a local holomorphic lift via the symmetrization map. In particular, the isometry $F_N$ is non-liftable.

\item We next establish the non-existence of isometric embeddings for $k=2$, see Theorem \ref{d=1}. This is achieved by first showing that any such isometry must intersect the royal variety. A holomorphic map from $\mathbb D$ into $\mathbb G_2$ whose boundary cluster set is contained in the distinguished boundary, extends holomorphically across $\partial\mathbb D$ and has rational components, see Proposition \ref{cluster_set}. This result plays a crucial role in proving that any isometry for $k=2$ cannot intersect the royal variety transversely. This obtains a contradiction and rules out the case $k=2$. 

\item We conclude by establishing the uniqueness of the non-royal isometric embedding for $k=1$. We first show that the image of such an isometry is contained in an irreducible projective algebraic curve $\Gamma$. However, this is a well-known phenomenon in literature \cite{Mok}, and follows from the rationality of the Bergman kernels of $\mathbb{D}$ and $\mathbb{G}_2$, together with the functional equation. 
The first major step in proving the uniqueness is to establish that the curve $\Gamma$ must be a projective line. The second major step is to establish that, after a possible reparametrization, the projective line $\Gamma$ corresponds to $S = -\sqrt{2}P$. This gives the non-royal isometry up to a rotation of the unit disk, thereby completing the classification.
\end{enumerate}

\medskip

The paper is organized as follows. In Section 2, we derive the functional equation satisfied by a holomorphic isometric embedding. In Section 3, we study the boundary cluster set and obtain restrictions on the isometric constant. Section 4 contains the royal and non-royal examples, together with a proof that they are not totally geodesic. In Section 5, we characterize the isometries for $k=3$, establish an obstruction to liftability for $k=1$ and $k=2$, and prove the nonexistence of isometries for $k=2$. Finally, in Section 6, we establish the uniqueness of the non-royal isometry for $k=1$, thereby completing the classification of holomorphic isometric embeddings from the unit disk into the symmetrized bidisk.


\section{Functional equation}

We shall begin by deriving the functional equation associated with the isometric embedding $F: (\mathbb{D}, kds^2)\rightarrow (\mathbb{G}_2, ds^2)$, which will play a central role in our subsequent analysis. Let $s=\lambda_1 + \lambda_2$ and $p = \lambda_1 \lambda_2$. It is easy to see that
\[
\mathbb{G}_2 = \left\{(s,p)\in \mathbb{C}^2: \max\left\{\left\lvert\frac{s \pm\sqrt{s^2-4p}}{2}\right\rvert\right\}<1\right\}.
\]
It follows from \cite{EZ} (see also \cite{MSZ}) that the Bergman kernel of the symmetrized bidisk is given by
\[
K_{\mathbb{G}_2}((s,p)) 
=
\frac{2 - \vert \lambda_1 + \lambda_2\vert^2 + 2 \vert\lambda_1\vert^2\vert\lambda_2\vert^2}{\pi^2 ((1-\vert\lambda_1\vert^2) (1-\vert\lambda_2\vert^2)(1 - \lambda_1\ov{\lambda_2})(1-\lambda_2\ov{\lambda_1}))^2}
=
\frac{2-\vert s\vert^2 + 2\vert p\vert^2}{\pi^2 \,\Gamma^2}    
\]
where
\begin{eqnarray*}
\Gamma
&=&
(1-\vert\lambda_1\vert^2)(1-\vert\lambda_2\vert^2)(1 -\lambda_1\ov{\lambda_2})(1-\lambda_2\ov{\lambda_1})
\\&=&
(1-\vert\lambda_1+\lambda_2\vert^2 + 2 \Re(\lambda_1\overline{\lambda_2}) + \vert\lambda_1\lambda_2\vert^2) (1-2 \Re(\lambda_1\overline{\lambda_2})+\vert\lambda_1\lambda_2\vert^2).
\end{eqnarray*}
Here
\[
\Re(\lambda_1\overline{\lambda_2}) = \Re\left( \left(\frac{s+\sqrt{s^2-4p}}{2}\right)\left(\overline{\frac{s-\sqrt{s^2-4p}}{2}}\right)\right)
=
\frac{\vert s\vert^2 - \vert s^2 - 4p\vert}{4}.
\]
So,
\begin{eqnarray*}
\Gamma &=& \left(1-\vert s\vert^2 + \frac{\vert s\vert^2 - \vert s^2 - 4p\vert}{2} + \vert p\vert^2\right) 
\left( 1 - \frac{\vert s\vert^2 - \vert s^2 - 4p\vert}{2} + \vert p\vert^2\right) 
\\&=& \frac{1}{4}
\left((2-\vert s\vert^2 + 2\vert p\vert^2)^2 - \vert s^2 - 4p\vert^2\right).
\end{eqnarray*}
Therefore,
\begin{equation}
K_{\mathbb{G}_2}((s,p)) = \frac{16}{\pi^2} \left(\frac{2-\vert s\vert ^2 + 2\vert p \vert^2}{((2-\vert s\vert^2 + 2\vert p\vert^2)^2 - \vert s^2 - 4p\vert^2)^2}\right).
\end{equation}

\medskip

\begin{rem}
For every automorphism $\varphi\in \operatorname{Aut}(\mathbb D)$, define $\Phi_\varphi:\mathbb G_2\to \mathbb G_2$
by
\[
\Phi_\varphi\left(\pi(z_1,z_2)\right)
=
\pi\left(\varphi(z_1),\varphi(z_2)\right).
\]
The maps $\Phi_\varphi$ give all the automorphisms of $\mathbb G_2$, and the inverse of $\Phi_{\phi}$ is $\Phi_{\varphi^{-1}}$. The automorphism group of $\mathbb{G}_2$ is not transitive. Thus, we may not always assume, without loss of generality, that the isometry $F$ vanishes at the origin. However, we can always assume, without loss of generality, that either $s_0 = 0$ or $p_0 = 0$.
This is well known; see \cite{TB, JP_aut} for details.
\end{rem}

Given a bounded domain $D \subset \mbb C^n$ and $\mu\in C(\ov{D})$ such that $\mu$ is positive on $D$, let $L^2_{\mu}(D)$ denote the space of all functions on $D$ that are square integrable with respect to $\mu dV$. Then, $A^2_{\mu}(D) = L^2_{\mu}(D) \cap \mathcal O(D)$ is a Hilbert space, and for any $z \in D$, the point evaluations $f \mapsto f(z) \in \mbb C$ are bounded (see \cite{PW2}, \cite{PW1}). So, there is a reproducing kernel $K_{D, \mu}(z, w)$, called the weighted Bergman kernel with respect to the weight $\mu$. Moreover, 
\[
K_{D,\mu}(z) = \sup\{\vert f(z) \vert^2 : f \in A^2_{\mu}(D), \Vert f \Vert_{A^2_\mu(D)} = 1 \}
\]
is a positive plurisubharmonic function on $D$. The weighted Bergman metric 
\[
ds^2_{D, \mu} = \sum_{i, j} \frac{\pa^2}{\pa z_i \pa \ov z_j} \log K_{D, \mu}(z) \; dz_i \otimes d \ov z_j
\]
is, therefore, well-defined on $D$. When $\mu \equiv 1$, we obtain the Bergman metric $ds^2_D$ on $D$. 

\medskip

On the unit disk $\mathbb{D}$, consider weights of the form $\mu_d(z) = K_{\mathbb{D}}^{-d}(z)$ for some integer $d \ge 0$. The Bergman kernel blows up at the boundary and hence $\mu_d$ can be extended continuously to the boundary $\pa \mathbb D$ by setting $\mu_d \equiv 0$ on $\pa \mathbb D$. The corresponding weighted Bergman kernel will be denoted by $K_{d}$ and the weighted Bergman metric by $ds^2_d$. It is easy to see that
\[
ds^2_d = (d+1) \, ds^2.
\]
Therefore, we can write $F:(\mathbb{D}, ds^2_{d}) \rightarrow (\mathbb{G}_2,ds^2)$, where $d=k-1$. In order to obtain the functional equation in full generality, we shall use the next theorem.

\begin{thm}[see \cite{JV}]\label{iso}
Let $D_1, D_2 \subset \mbb C^n$ be bounded domains and $\mu_1, \mu_2$ are weights as above. Let $f : (D_1, ds^2_{D_1, \mu_1}) \ra (D_2, ds^2_{D_2, \mu_2})$
be a holomorphic isometry embedding. Then, there exist functions $s_0 \in A^2_{\mu_1}(D_1)$ and $r_0 \in A^2_{\mu_2}(D_2)$ such that
\[
s_0(z) \; \ov{s_0(w)} \; K_{D_2, \mu_2}(f(z), f(w)) = r_0(f(z)) \; \ov{r_0(f(w))} \; K_{D_1, \mu_1}(z, w)
\]
for all $z, w \in D_1$. 
\end{thm}

It follows from the proof of Theorem \ref{iso} that the functions $s_0$ and $r_0$ can be chosen as follows: For a fixed $p\in D_1$, choose orthonormal bases $\{s_i\}_{i=0}^{\infty}$ and $\{r_i\}_{i=0}^{\infty}$ of $A^2_{\mu_1}(D_1)$ and $A^2_{\mu_2}(D_2)$ respectively so that $s_i(p)=0$, $r_i(f(p))=0$ for $ i\geq 1$, and $s_0(p)\neq 0$, $r_0(f(p))\neq 0$. 

\medskip

Note that $\mu_d(z) = \pi^d (1-\vert z\vert^2)^{2d}$. An orthonormal basis of $A^2_{d}(\mathbb{D})$ is given by (see \cite{JV})
\[
\left\{\sqrt{\frac{1}{\pi^{d+1}}\frac{(m+2d+1)!}{m! (2d)!}} \,\,z^m\right\}_{m=0}^{\infty}.
\]
The weighted Bergman kernel $K_d$ is given by (see \cite{JV})
\[
K_d(z,w) = \frac{(2d+1)}{\pi^{d+1}}\frac{1}{(1-z\ov{w})^{2d+2}}.
\]
An orthonormal basis of $A^2(\mathbb{G}_2)$ can be defined using Schur polynomials \cite{MSZ}. A partition $\mathbf{p}$ is a tuple $\mathbf{p}:= (p_1,p_2)$ of non-negative integers such that $p_1\geq  p_2$. Let $[[2]]$ denote the set of all partitions of the form $\mathbf{p} = \mathbf{q} + \delta$ where $\mathbf{q}$ is a partition and $\delta = (1,0)$. For $\mathbf{p}\in[[2]]$, the Schur polynomials are
\[
S_{\mathbf{p}}(z)  =   \frac{z_1^{p_1}z_2^{p_2} - z_2^{p_1}z_1^{p_2}}{z_1-z_2}, \quad\quad
z = (z_1,z_2)\in\mathbb{D}_2.
\]
The Schur polynomials are indeed polynomials and they are symmetric. Therefore, there exist a holomorphic function $f_{\mathbf{p}}:\mathbb{G}_2\rightarrow\mathbb{C}$ such that $f_{\mathbf{p}}\circ \pi = S_{\mathbf{p}}$. Then an orthonormal basis of $A^2(\mathbb{G}_2)$ is given by
\[
\left\{\sqrt{\frac{(p_1+1)(p_2+1)}{\pi^2}} \,\,f_{\mathbf{p}}: \mathbf{p}\in[[2]]\right\}.
\]
Note that $f_{(1,0)}$ is a constant equal to $\sqrt{2}/\pi$ and other functions in the orthonormal basis vanish at $(0,0)$. 

\subsection{Functional equation when $\mathbf{F(0) = (0,0)}$.}

In Theorem \ref{iso}, taking 
\[
s_0 = \sqrt{\frac{2d+1}{\pi^{d+1}}}
\quad\text{and}\quad
r_0 = f_{(1,0)} = \frac{\sqrt{2}}{\pi}
\]
we obtain
\begin{eqnarray*}
K_{\mathbb{G}_2}(F(z),F(w)) 
&=&
\frac{2}{\pi^2}\,\frac{\pi^{d+1}}{(2d+1)}\,K_d(z,w) 
=
\frac{2}{\pi^2} \frac{1}{(1-z\ov{w})^{2d + 2}}
\\ &=&
\frac{2}{\pi^2} \frac{1}{(1-z\ov{w})^{2k}}
\end{eqnarray*}
for $z,w\in\mathbb{D}$. Along the diagonal, we can write
\begin{equation}\label{main}
K_{\mathbb{G}_2}(F(z)) = \frac{2}{\pi^2}\frac{1}{(1-\vert z\vert^2)^{2k}}.
\end{equation}

\subsection{Functional equation in the general case}

Let $F(0) = (s,p)\in\mathbb{G}_2$. 

\medskip

Let $\{r_i\}_{i=0}^{\infty}$ be an orthonormal basis of $A^2(\mathbb{G}_2)$ such that $r_0((s,p)) \neq 0$ and $r_i((s,p))= 0$ for all $i\geq 1$. Taking $s_0 = \sqrt{2d+1}/\sqrt{\pi^{d+1}}$, it follows from Theorem \ref{iso} that
\begin{eqnarray*}
K_{\mathbb{G}_2}(F(z), F(w)) 
&=&
r_0(F(z)) \,\overline{r_0(F(w))}\,\frac{\pi^{d+1}}{2d+1}\, K_{\mathbb{D},d}(z)  
\\&=&
r_0(F(z)) \,\overline{r_0(F(w))} \,\frac{1}{(1- z\overline{w})^{2k}}.
\end{eqnarray*}
We therefore need to find out the function $r_0$. Consider the evaluational functional
\[
E_{(s,p)}: A^2(\mathbb{G}_2)\rightarrow \mathbb{C}
\quad 
\text{given by } f\mapsto f((s,p)).
\]
The orthogonal complement of the kernel of $E_{(s,p)}$ is one-dimensional. Let $A=\{ f\in A^2(\mathbb{G}_2): f((s,p)) = 0\}$ denote the kernel of $E_{(s,p)}$. So, $r_0$ is the function in the orthogonal complement of $A$ with norm $1$. For every $f\in A^2(\mathbb{G}_2)$
\[
0 = \langle r_0, f-f((s,p))\rangle = \langle r_0, f\rangle - \langle r_0, f((s,p))\rangle.
\]
So,
\begin{eqnarray*}
r_0((z_1,z_2)) 
&=&
\langle r_0, K_{\mathbb{G}_2}(\cdot, (z_1,z_2))\rangle 
\\&=&
\langle r_0, K_{\mathbb{G}_2}((s,p), (z_1,z_2)) \rangle 
=
K_{\mathbb{G}_2}((z_1,z_2),(s,p)) \int_{\zeta\in \mathbb{G}} r_0(\zeta) dV(\zeta).
\end{eqnarray*}
Therefore, we see that
\[
r_0((z_1,z_2)) = c \,K_{\mathbb{G}_2}((z_1,z_2),(s,p))
\]
for some constant $c$. Since $r_0$ has norm $1$, we see that
\[
1 = \langle r_0, r_0\rangle = \vert c\vert^2 \langle K_{\mathbb{G}_2}(\cdot,(s,p)), K_{\mathbb{G}_2}(\cdot,(s,p)) \rangle = \vert c\vert^2 K_{\mathbb{G}_2}((s,p),(s,p)).
\]
Therefore, 
\[
r_0 = \frac{K_{\mathbb{G}_2}(\cdot,(s,p))}{\sqrt{K_{\mathbb{G}_2}((s,p),(s,p))}}.
\]
This gives us the functional equation
\[
K_{\mathbb{G}_2}(F(z), F(w)) 
=
\frac{1}{K_{\mathbb{G}_2}((s,p),(s,p))} K_{\mathbb{G}_2}(F(z),(s,p)) K_{\mathbb{G}_2}((s,p), F(w)) \frac{1}{(1-z\overline{w})^{2k}}.
\]
We can also write it as
\begin{equation}
\frac{K_{\mathbb{G}_2}(F(z), F(w))}{K_{\mathbb{G}_2}(F(z),(s,p)) \, K_{\mathbb{G}_2}((s,p), F(w))} 
=
\frac{1}{K_{\mathbb{G}_2}((s,p),(s,p))}\frac{1}{(1-z\overline{w})^{2k}}.
\end{equation}
This can be written along the diagonal as
\begin{equation}\label{feqn}
\frac{K_{\mathbb{G}_2}(F(z))}{\left\lvert K_{\mathbb{G}_2}(F(z),(s,p))\right\rvert^2} 
=
\frac{1}{K_{\mathbb{G}_2}((s,p))}\frac{1}{(1-\vert z\vert^2)^{2k}}.
\end{equation}
It is easy to see that the function $F$ must be a proper holomorphic map. 

\section{Isometric constant and boundary behavior}

Let $\partial_s\mathbb{G}_2$ denote the distinguished boundary of $\mathbb{G}_2$. It is known that
\[
\partial_s\mathbb{G}_2 = \pi(\partial \mathbb{D} \times \partial \mathbb{D}).
\]
The set $\mathcal{R} = \{\pi(z,z): z\in \overline{\mathbb{D}}\}$ is called the \textit{royal variety} of $\mathbb{G}_2$.
Divide $\partial\mathbb{G}_2$ into three disjoint subsets
\[
\mathcal{B}_1:= \partial \mathbb{G}_2\setminus\partial_s \mathbb{G}_2 = \{\pi(\lambda_1,\lambda_2): (\lambda_1,\lambda_2)\in(\mathbb{D}\times\partial\mathbb{D}) \cup (\partial\mathbb{D}\times \mathbb{D})\}
\]
\[
\mathcal{B}_2:=\partial_s\mathbb{G}_2\setminus\mathcal{R} = \{\pi(\lambda_1,\lambda_2) : \vert \lambda_1 \vert = \vert \lambda_2\vert =1, \lambda_1\neq \lambda_2\}.
\]
\[
\mathcal{B}_3:=\partial \mathbb{G}_2\cap \mathcal{R} = \{(2\lambda,\lambda^2): \vert\lambda\vert =1\}
\]
Let $F:(\mathbb D,k\,ds^2)\to(\mathbb G_2,ds^2)$
be a holomorphic isometric embedding. In this section,
we show that the isometric constant $k$ can take only
three possible values, namely $1$, $2$, and $3$.
Moreover, these values correspond to three distinct
types of boundary behavior.

\medskip

Let $A\approx B$ means that $A/B$ is bounded
above and below by positive constants.

\begin{prop}\label{approx}
Let
\[
F:(\mathbb D,k\,ds^2)
\rightarrow(\mathbb G_2,ds^2)
\]
be a holomorphic isometric embedding. Then
\begin{equation}
K_{\mathbb G_2}(F(z))
\approx
\frac{1}{(1-|z|^2)^{2k}},
\qquad z\in\mathbb D.
\end{equation}
\end{prop}

\begin{proof}
We can assume, without loss of generality, that $F(0) = (s_0,0)$.
The functional equation gives
\[
\frac{K_{\mathbb G_2}(F(z))}
     {|K_{\mathbb G_2}(F(z),(s_0,0))|^2}
=
\frac{1}{K_{\mathbb G_2}((s_0,0))}
\frac{1}{(1-|z|^2)^{2k}}.
\]
For $(s_1,p_1), (s_2,p_2)\in\mathbb{G}_2$, by polarizing the formula of the Bergman kernel of $\mathbb{G}_2$, we have
\[
K_{\mathbb{G}_2}((s_1,p_1),(s_2,p_2)) 
= 
\frac{16}{\pi^2}
\left(\frac{2 - s_1 \overline{s_2} + 2 p_1 \ov{p_2}}{((2 - s_1 \ov{s_2} + 2 p_1 \ov{p_2})^2 - (s_1^2 - 4 p_1) (\ov{s_2^2 - 4 p_2}))^2}
\right).
\]
Let $F = (S,P)$. Then,
\[
K_{\mathbb G_2}(F(z), (s_0,0))
=
\frac{16}{\pi^2}
\frac{2-\overline {s_0}\,S(z)}
{\left((2-\overline {s_0}\,S(z))^2 - \ov{s_0}^2(S(z)^2 - 4 P(z))\right)^2}
=
\frac{1}{\pi^2}
\frac{2-\overline {s_0}\,S(z)}
{\left( 1 - \overline {s_0}\,S(z) + \ov{s_0}^2 P(z)\right)^2}.
\]
For a fixed $z\in\mathbb D$, choose $\lambda_1,\lambda_2\in\mathbb D$ such that
$S(z)=\lambda_1+\lambda_2$ and $P(z)=\lambda_1\lambda_2$.
Since
\[
1-\overline{s_0}S(z)+\overline{s_0}^{\,2}P(z)
=
(1-\overline{s_0}\lambda_1)
(1-\overline{s_0}\lambda_2),
\]
we have
\[
(1-|s_0|)^2
\leq
|1-\overline{s_0}S(z)+\overline{s_0}^{\,2}P(z)|
\leq
(1+|s_0|)^2.
\]
Moreover,
\[
2(1-|s_0|)
\leq |2-\overline{s_0}S(z)|
\leq 2(1+|s_0|).
\]
It follows that
\[
\frac{2(1-|s_0|)}{\pi^2(1+|s_0|)^4}
\leq
|K_{\mathbb G_2}(F(z),(s_0,0))|
\leq
\frac{2(1+|s_0|)}{\pi^2(1-|s_0|)^4}.
\]
Since $z\in\mathbb{D}$ was arbitrary, we have
\[
K_{\mathbb G_2}(F(z)) \approx \frac{1}{(1-|z|^2)^{2k}}.
\]
This proves the lemma.
\end{proof}

The boundary cluster set of the holomorphic isometric embedding $F$ is given by
\[
C(F)=
\bigcap_{0<r<1}
\overline{F\bigl(\{z\in\mathbb D:|z|>r\}\bigr)}.
\]

\begin{prop}\label{cluster_set}
Let $F:=(S,P):\mathbb D
\rightarrow\mathbb G_2$
be a holomorphic function. If 
\[
C(F)\subset\partial_s\mathbb G_2,
\]
then $F$ extends holomorphically across
every point of $\partial\mathbb D$. Moreover, $S$ and $P$ are rational functions.
\end{prop}

\begin{proof}
Since $C(F)\subset\partial_s\mathbb G_2$, we have $\vert P(z)\vert \rightarrow 1$ as $\vert z\vert \rightarrow 1$. Therefore, $P:\mathbb{D}\rightarrow\mathbb{D}$ is a proper holomorphic map and hence a finite Blaschke product. The proof is inspired from the description of rational-inner functions for the symmetrized bidisk, due to Agler-Lykova-Young \cite{Agler}. Let $N = \deg P$.
We can write
\[
P(z)=\tau \frac{\widetilde D(z)}{D(z)},
\qquad |\tau|=1,
\]
where
\[
D(z)=\prod_{j=1}^{N}(1-\overline{\alpha_j}z),
\qquad
\widetilde D(z)=z^{N}\overline{D(1/\overline z)},
\qquad
\alpha_j\in\mathbb{D}.
\]
Since $S$ is bounded and holomorphic, its radial boundary values exist almost everywhere. As all boundary cluster values belong to $\partial_s\mathbb G_2$, we have
\[ 
S=\overline S\,P
\qquad \text{a.e. on } \partial\mathbb D.
\]
The idea is to prove that $S$ is also a rational function of the form $E/D$ for some polynomial $E$. Put $Q:=DS$. The function $Q$ is bounded and holomorphic on $\mathbb{D}$.
For $\zeta\in\partial\mathbb D$, we have
\[
P(\zeta)=\tau \zeta^N\frac{\overline{D(1/\ov \zeta)}}{D(\zeta)}
=
\tau \zeta^N\frac{\overline{D(\zeta)}}{D(\zeta)}.
\]
Therefore,
\[
Q(\zeta)=
D(\zeta) S(\zeta)
=
D(\zeta) \ov{S(\zeta}) P(\zeta)
=
\tau \zeta^N \overline{Q(\zeta)}
\qquad \text{a.e. on } \partial\mathbb D.
\]
We write the power series representation of $Q$ as
\[
Q(z) = \sum_{k=0}^{\infty} q_k z^k.
\]
Comparing the Fourier coefficients on both sides of
\[
Q(\zeta) = \tau \zeta^N \overline{Q(\zeta)},
\]
we conclude that $Q$ is a polynomial of degree at most $N$. Therefore, $S = Q/D$ is rational.
\end{proof}

\begin{prop}\label{isoc}
Let
\[
F:(\mathbb D,k\,ds_{\mathbb D}^2)
\rightarrow
(\mathbb G_2,ds_{\mathbb G_2}^2)
\]
be a holomorphic isometric embedding, where
$k$ is a positive integer. Then $C(F)\subset\partial\mathbb G_2$ and $k\in\{1,2,3\}$.
Moreover,
\[
\begin{aligned}
k=1
&\iff C(F)\cap \mathcal B_1\ne\varnothing,\\
k=2
&\iff C(F)\subset \mathcal B_2,\\
k=3
&\iff C(F)\subset \mathcal B_3.
\end{aligned}
\]
\end{prop}

\begin{proof}
It follows from the functional equation \ref{feqn} of the isometry that $F$ is a proper holormorphic map. Therefore, $C(F) \subset \partial\mathbb{G}_2$.
We first establish an estimate for the roots
of the polynomial
\[
\lambda^2-S(z)\lambda+P(z)=0,
\]
where $F=(S,P)$. For each $\omega\in\partial\mathbb D$, consider the maps
\[
\Phi_\omega(s,p)
=\frac{2\omega p-s}{2-\omega s}.
\]
It is known \cite{Agler_0} that $\Phi_\omega$ maps $\mathbb G_2$ holomorphically into $\mathbb D$. Since $\partial\mathbb D$ is compact, we have $m:=\max_{|\omega|=1}|\Phi_\omega(F(0))|<1$.
The Schwarz--Pick lemma applied to
$\Phi_\omega\circ F$ gives
\[
\left\lvert\frac{(\Phi_\omega\circ F)(z) - (\Phi_\omega\circ F)(0)}{1 - \ov{(\Phi_\omega\circ F)(0)} \,(\Phi_\omega\circ F) (z)} \right\rvert \leq \vert z\vert
\]
which by applying triangle inequalities gives
\[
\vert (\Phi_\omega\circ F)(z)\vert 
\leq
\frac{\vert z\vert + \vert (\Phi_\omega\circ F)(0)\vert}{1 + \vert z\vert\, \vert (\Phi_\omega\circ F)(0)\vert}
\leq
\frac{\vert z\vert + m}{1+m \vert z\vert}.
\]
Therefore,
\[
1-|\Phi_\omega(F(z))|
\geq
\frac{1-m}{1+m \vert z\vert}(1-|z|)
\geq
\frac{1-m}{1+m}(1-|z|)
\]
uniformly in $\omega\in\partial\mathbb D$. For a fixed $z\in\mathbb{D}$, if $\lambda_1,\lambda_2\in\mathbb{D}$ are the solutions of $\lambda^2 - S(z) \lambda + P(z) = 0$, put $\omega = \ov{\alpha_1}$, where $\lambda_1 = \vert \lambda_1\vert \alpha_1$. Now,
\[
\Phi_{\om}(F(z)) 
=
\frac{2 \overline{\alpha_1}\lambda_1\lambda_2 - \lambda_1 - \lambda_2}{2 - \ov{\alpha_1}(\lambda_1 + \lambda_2)}
=
\alpha_1 \left(
\frac{2 \vert \lambda_1\vert\overline{\alpha_1} \lambda_2 - \vert \lambda_1\vert - \ov{\alpha_1}\lambda_2}{2 - \vert \lambda_1\vert - \ov{\alpha_1}\lambda_2}
\right).
\]
Note that when $\vert\lambda_1\vert \rightarrow 1$, the right hand side approaches $-\alpha_1$. Therefore, consider
\begin{eqnarray*}
\left \lvert \Phi_{\om}(F(z)) + \alpha_1 \right\rvert
&=&
\left\lvert
2\alpha_1 
\left(
\frac{\vert \lambda_1\vert \ov{\alpha_1} \lambda_2 - \ov{\alpha_1} \lambda_2 + 1 - \vert \lambda_1\vert}{2 - \vert \lambda_1\vert - \ov{\alpha_1}\lambda_2}
\right)
\right\rvert
\\
&=&
2
\left\lvert 
\frac{(1 - \vert \lambda_1\vert ) ( 1 - \overline{\alpha_1} \lambda_2)}{(1 - \vert \lambda_1\vert) + (1 - \ov{\alpha_1}\lambda_2)}
\right\rvert
\\&\leq &
2
( 1 - \vert \lambda_1\vert ),
\end{eqnarray*}
where the last inequality holds because $\Re ( 1 - \ov{\alpha_1} \lambda_2) > 0$. Therefore,
\begin{eqnarray*}
1 - \vert \lambda_1 \vert^2
\geq 
1 - \vert \lambda_1\vert
&\geq&
\frac{1}{2} \left \lvert \Phi_{\om}(F(z)) + \alpha_1 \right\rvert
\\&\geq&
\frac{1}{2} \left( \vert \alpha_1 \vert - \vert \Phi_{\omega}(F(z)\vert\right)
=
\frac{1}{2} \left( 1 - \vert \Phi_{\omega}(F(z)\vert\right)
\\&\geq&
\frac{1}{2} \frac{1 - m}{1 + m} (1 - \vert z\vert)
\\& \geq &
\frac{1}{4} \frac{1 - m}{1 + m}\, (1 - \vert z\vert^2).
\end{eqnarray*}
Similar calculations will hold for $\lambda_2$. Therefore, we obtain
\begin{equation}
1 - \vert \lambda_i \vert^2 
\geq 
\frac{1}{4} \frac{1 - m}{1 + m} (1 - \vert z\vert^2)
\quad
\text{ for } i = 1,2,
\end{equation}
where $z\in\mathbb{D}$ satisfies $F(z) = (\lambda_1 + \lambda_2, \lambda_1 \lambda_2).$ Moreover, 
\begin{equation}\label{kernel}
K_{\mathbb{G}_2}(F(z)) 
=
\frac{(1-\vert\lambda_1\vert^2) (1-\vert\lambda_2\vert^2) + \vert 1 - \lambda_1\ov{\lambda_2}\vert^2}
{\pi^2 ((1-\vert\lambda_1\vert^2) (1-\vert\lambda_2\vert^2)\vert 1 - \lambda_1\ov{\lambda_2}\vert^2)^2}.
\end{equation}
Let
\[
A=(1-|\lambda_1|^2)(1-|\lambda_2|^2),
\qquad
B=|1-\lambda_1\overline \lambda_2|^2.
\]
Since $B=A+|\lambda_1-\lambda_2|^2\ge A$, we have
\[
K_{\mathbb G_2}(F(z))
=
\frac{A+B}{\pi^2A^2B^2}
\le 
\frac{2}{\pi^2A^3}
\le
\frac{C_0}{(1-|z|^2)^6},
\]
for some constant $C_0>0$, which is independent of $z\in\mathbb{D}$. Combining this with Proposition \ref{approx}, we obtain $2k \leq 6$, that is, $k\in\{1, 2, 3\}$.

\medskip

Next, we examine the possible locations of the cluster values. If $C(F) \subset \mathcal{B}_2 \cup \mathcal{B}_3 = \partial_s\mathbb{G}_2$, then $F$ extends holomorphically across $\zeta\in\partial\mathbb{D}$. Therefore, if $F(\zeta)\in \mathcal B_2$, we have
\[
A(r\zeta)\approx(1-r^2)^2,
\qquad B(r\zeta)\approx 1,
\]
giving
\[
K_{\mathbb G_2}(F(r\zeta))
\approx (1-r^2)^{-4}.
\]
Thus, $k=2$. If $F(\zeta)\in \mathcal B_3$, we have
\[
A(r\zeta)\approx(1-r^2)^2,
\quad B(r\zeta)\approx (1-r^2)^2
\quad
\text{and}
\quad
(A + B)(r\zeta) \approx (1-r^2)^2
\]
giving
\[
K_{\mathbb G_2}(F(r\zeta))
\approx(1-r^2)^{-6}.
\]
Here, we obtain $k=3$.

\medskip

Now, suppose that $\zeta\in C(F)\cap \mathcal B_1$. Choose a sequence
$z_n\to\partial\mathbb D$ such that
$F(z_n)\to \zeta$. After interchanging the
roots if necessary, write $F(z_n)=\pi(a_n,b_n)$, where
\[
a_n\to\lambda,\qquad b_n\to\mu,
\qquad |\lambda|=1,\quad |\mu|<1.
\]
Then,
\[
A(z_n) \geq C_1 (1 - \vert z_n\vert^2)
\quad\text{and}\quad
B(z_n) \geq C_2
\]
for some constants $C_1, C_2>0$. Therefore, we have
\[
K_{\mathbb G_2}(F(z_n))
\le \frac{2}{\pi^2\,C_1^2 \,C_2}\frac{1}{(1-|z_n|^2)^2}.
\]
Combining this with Proposition \ref{approx}, we have $k = 1$. Conversely, if $k=1$ and $C(F)\cap \mathcal{B}_1 = \varnothing$, it would imply that $C(F) \subset \mathcal{B}_2 \cup \mathcal{B}_3$. But this would force $k=2$ or $k=3$, which is a contradiction. Therefore, 
\[
k=1 \quad \text{if and only if} \quad C(F)\cap \mathcal{B}_1 \neq \varnothing.
\]
Thus, $k=2$ or $k=3$ implies that $C(F) \cap \mathcal{B}_1 = \varnothing$. This completely proves the other implications
\[
k = 2 \iff C(F)\subset\mathcal{B}_2
\quad
\text{and}
\quad
k = 3 \iff C(F) \subset\mathcal{B}_3.
\]
This finishes the proof.
\end{proof}


\section{A non-trivial isometry}

In this section, we present two examples of holomorphic isometric embeddings $F:(\mathbb{D}, kds^2) \rightarrow (\mathbb{G}_2,ds^2)$ for $k=1$ and $k=3$ respectively, and show that neither is totally geodesic. We will see in the next section that there does not exist any holomorphic isometric embedding for $k=2$.

\medskip

Our first example is obtained by composing the diagonal map $z\mapsto (z,z)$ with the symmetrization map $\pi$. It is known \cite{Ng} that $G:(\mathbb{D}, 2 ds^2)\rightarrow (\mathbb{D}\times \mathbb{D}, ds^2)$ defined by
\[
G : z\mapsto (z,z)
\]
is a holomorphic isometric embedding. The isometry $G$ is totally geodesic. However, as we shall see, $\pi\circ G$ is not totally geodesic.

\begin{exam}[Royal isometry]
Define $F:=(S,P):(\mathbb{D},3ds^2)\rightarrow (\mathbb{G}_2,ds^2)$ by
\begin{equation}
F(z) = \pi(z,z) = (2z, z^2).
\end{equation}
The function $F$ is clearly well-defined. Note that, $S^2 \equiv 4 P$. Therefore,
\begin{eqnarray*}
K_{\mathbb{G}_2}(F(z)) 
&=&
\frac{16}{\pi^2} \left(
\frac{2 - 4\vert z\vert^2 + 2 \vert z\vert^4}{(2 - 4\vert z\vert^2 + 2\vert z\vert^4)^4}
\right)
\\&=&
\frac{1}{\pi^2} 
\frac{1}{(1 - 2\vert z\vert^2 + \vert z\vert^4)^3}
\\&=&
\frac{1}{\pi^2} 
\frac{1}{(1 - \vert z\vert^2)^6}.
\end{eqnarray*}
This proves that $F:=(S,P):(\mathbb{D},3ds^2)\rightarrow (\mathbb{G}_2,ds^2)$ is a holomorphic isometric embedding.
\end{exam}

\medskip

Our next example has this interesting property that it is not liftable; we shall prove this fact in Corollary \ref{nolift}.

\begin{exam}[Non-royal isometry]\label{nonroyal}
Define $F:=(S,P):(\mathbb{D},ds^2)\rightarrow (\mathbb{G}_2,ds^2)$ by
\begin{equation}
F(z) := (S,P)= \left(\sqrt{2}(\sqrt{z+1} - 1), 1-\sqrt{z+1}\right)
\end{equation}
Here, $F(0) = (0,0)$. For $z\in \mathbb{D}$, we have $z+1\in \mathbb{C}\setminus(-\infty,0]$. Therefore, $\sqrt{z+1}$ is well-defined. Note that, 
\[
S = -\sqrt{2}P. 
\]
We first prove that $\vert P\vert < 1$. For $z\in\mathbb{D}$,
\[
1> \vert z\vert =
\vert (1-P(z))^2-1\vert = 
\vert P(z) (2-P(z))\vert.
\]
Now, $2-P(z) = 1 + \sqrt{z+1}$. For $z\in \mathbb{D}$, we have $\arg(z+1)\in (-\pi/2,\pi/2)$. Thus, $\Re(\sqrt{z+1}) > 0$. This implies that $\vert P(z)\vert = \vert 1 - \sqrt{z+1}\vert\leq \vert 1+ \sqrt{z+1}\vert$. Therefore, $\vert P(z)\vert^2<1$. This proves that $\vert P(z)\vert <1$ for $z\in\mathbb{D}.$ 

\medskip

Let $A(z) = \sqrt{z+1}$ for $z\in\mathbb{D}$. Then,
\begin{eqnarray*}
\vert S-\overline{S}P\vert^2 - (1-\vert P\vert^2)^2
&=&
\vert -\sqrt{2} P + \sqrt{2}\vert P\vert^2\vert^2 - (1-\vert P\vert^2)^2
\\&=&
2 \vert P\vert^2 \vert 1- P\vert^2 - (1-\vert P\vert^2)^2
\\&=&
2 (1+\vert A\vert^2 - 2\Re(A)) \vert A\vert^2 - (2\Re(A)-\vert A\vert^2)^2
\\&=&
2\vert A\vert^2 + \vert A\vert^4 - 4 (\Re(A))^2 
\\&=&
\vert A\vert^4 - 2 \Re(A^2) + 1 - 1
\\&=&
\vert A^2 - 1\vert^2 - 1 
=
\vert z\vert^2 - 1 < 0.
\end{eqnarray*}
Thus, $\vert S- \overline{S}P\vert < 1 - \vert P\vert^2$. Therefore, $F$ is well-defined.

\medskip

Now,
\begin{eqnarray*}
K_{\mathbb{G}_2}(F(z))
&=&
K_{\mathbb{G}_2}(S(z),P(z))
=
\frac{16}{\pi^2}\left(
\frac{2 - \vert S(z)\vert^2 + 2\vert P(z)\vert^2}{((2-\vert S(z)\vert^2 + 2\vert P(z)\vert^2)^2-\vert S(z)^2 - 4P(z)\vert^2)^2}
\right)
\\&=&
\frac{16}{\pi^2}\left(
\frac{2}{(4-\vert 2 P(z)^2 - 4P(z)\vert^2)^2}
\right)
\\&=&
\frac{2}{\pi^2}\left(
\frac{1}{(1-\vert P(z)\vert^2 \vert 2 - P(z)\vert^2)^2}
\right)
\\&=&
\frac{2}{\pi^2}\left(
\frac{1}{(1-\vert 1 - A(z)\vert^2 \vert 1 + A(z)\vert^2)^2}
\right)
=
\frac{2}{\pi^2}\left(
\frac{1}{(1-\vert 1 - (A(z))^2\vert^2 )^2}
\right)
\\&=&
\frac{2}{\pi^2}\left(
\frac{1}{(1-\vert z\vert^2 )^2}
\right).
\end{eqnarray*}
This proves that $F:(\mathbb{D},ds^2)\rightarrow (\mathbb{G}_2, ds^2)$ is a holomorphic isometric embedding. Note that $S^2-4P = 2z =0$ only at $z=0$. Thus, $F(\mathbb{D})$ intersects the royal variety exactly at one point $(0,0)$.
We can compose $F$ with automorphisms to get new isometries. Their images will also intersect the royal variety at exactly one point. 

\medskip

We also describe the boundary cluster set of $F$. Since
$z\mapsto\sqrt{z+1}$ extends continuously to
$\overline{\mathbb D}$, the map $F$ also extends
continuously to $\overline{\mathbb D}$. Consequently,
\[
C(F)=F(\partial\mathbb D)
=\left\{
\left(
\sqrt{2}(\sqrt{\zeta + 1}-1),
1-\sqrt{\zeta + 1}
\right):|\zeta|=1
\right\}.
\]
Let $\zeta=e^{i\theta}$, where
$-\pi<\theta<\pi$. Then
\[
1 + \zeta 
=
1 + \cos{\theta} + i \sin\theta
=
2 \cos^2\left(\frac{\theta}{2}\right) 
+
i \,
2 \cos \left(\frac{\theta}{2}\right) \sin\left(\frac{\theta}{2}\right)
=
2 \cos\left(\frac{\theta}{2}\right) 
e^{i\theta/2}.
\]
Therefore, $\sqrt{\zeta + 1} = \sqrt{2 \cos(\theta/2)}\,e^{i\theta/4}$.
The roots of $\lambda^2-S(\zeta)\lambda+P(\zeta)=0$ are
\[
\lambda_\pm
=\frac{\sqrt{\zeta+1}-1\pm\sqrt{\zeta}}{\sqrt{2}}.
\]
Thus,
\[
\lambda_+
=\frac{e^{i\theta/4}}{\sqrt{2}}
\left(\sqrt{2 \cos(\theta/2)}+2i\sin(\theta/4)\right)
\quad
\text{and}
\quad
\lambda_-
=-\frac{\sqrt{2}e^{i\theta/4}}
{2\cos(\theta/4)+\sqrt{2 \cos(\theta/2)}}.
\]
Since $2 \cos(\theta/2) = 1 - 2 \sin^2(\theta/4)$, we have $\vert \lambda_+\vert^2 = 1$. And $2 \cos(\theta/4) + \sqrt{2 \cos(\theta/2)} \geq 2 \cos(\theta/4) > 2\cos (\pi/4) = \sqrt{2}$, which implies that $\vert \lambda_-\vert < 1.$
Therefore,
\[
F(e^{i\theta})\in \mathcal B_1
\quad
\text{for }
-\pi < \theta < \pi.
\]
However, 
\[
F(e^{i\pi}) = F(-1)=(-\sqrt{2},1) \in\mathcal{B}_2.
\]
Therefore,
\begin{equation}
C(F)\cap \mathcal B_2=\{(-\sqrt{2},1)\},
\quad\text{and}\quad
C(F)\setminus\{(-\sqrt{2},1)\}\subset \mathcal B_1.  
\end{equation}
\end{exam}

\subsection{Isometries are not totally geodesic} The Bergman metric of the unit disk is
\[
ds^2_{\mathbb{D}} = \frac{2}{(1-\vert z\vert^2)^2} dz \otimes d\ov{z}.
\]
Since all the first order partial derivatives of the coefficient of $ds^2_{\mathbb{D}}$ vanish at $0$, the Christoffel symbols of $ds^2_{\mathbb{D}}$ also vanish at $0$. Since the weighted Bergman metric $ds^2_{d}$ on the unit disk with respect to the weight $\mu_d$ is $(d+1) g_{\mathbb{D}}$, this shows that the Christoffel symbols of $ds^2_{d}$ must also vanish at $0$.

\medskip

The Christoffel symbols in holomorphic coordinates for the Bergman metric $ds^2$ on the symmetrized bidisk at $(0,0)$ can be read from \cite[Proposition 9]{CY} (see also \cite{Try}), and we have
\[
\Gamma^1_{11} = 0, 
\quad
\Gamma^2_{11} = \frac{-4}{5},
\quad
\Gamma^1_{12} = \Gamma^1_{21} = 0,
\quad
\Gamma^1_{22} = 0,
\quad
\Gamma^2_{21} = \Gamma^2_{12} = 0,
\quad
\Gamma^2_{22} = 0
\]

\medskip

Let $F:=(f_1,f_2): (\mathbb{D}, ds^2_{d})\rightarrow (\mathbb{G}_2,ds^2)$ be an isometric embedding. The isometry $F$ is said to be \textit{totally geodesic} when $(F(\mathbb{D}), F_*(ds^2_d))$ is a totally geodesic submanifold of $(\mathbb{G}_2, ds^2)$, that is, every $F_*(ds^2_d)$-geodesic in $F(\mathbb{D})$ is also a $ds^2$-geodesic in $\mathbb{G}_2$. This happens if and only if the second fundamental form $\mathrm{II}$ of $F(\mathbb{D})$ vanishes identically. If $X,Y \in \Gamma(T_{\mathbb{C}}(F(\mathbb{D})))$ are extended arbitrarily to complex vector fields on a neighborhood of $F(\mathbb{D})$ in $\mathbb{G}_2$, the Gauss Formula (see \cite[Theorem 8.2]{Lee_Riemannian} and \cite{Lee_complex}) states
\[
\mathrm{II}(X,Y) 
=
\nabla^{\mathbb{G}_2}_{X} Y - \nabla^{F(\mathbb{D})}_X Y.
\]
Let $X = Y = F_*\left(\frac{\partial}{\partial z}\right)$. Now,
\[
\nabla^{F(\mathbb{D})}_X \,Y = F_*\left(\nabla^{\mathbb{D}}_{(F^{-1})_*(X)} (F^{-1})_*(Y)\right) 
=
F_*\left(\nabla^{\mathbb{D}}_{\frac{\partial}{\partial z}}\, \frac{\partial}{\partial z}\right)
\]
vanishes at $0$ because all the Christoffel symbols of $ds^2_d$ vanish at $0$. Let $(s,p)$ denote the coordinates of $\mathbb{G}_2$. We have
\[
X = Y = F_*\left(\frac{\partial}{\partial z}\right) 
=
\left(\frac{\partial f_1}{\partial z} \circ F^{-1}\right) \frac{\partial}{\partial s} 
+
\left(\frac{\partial f_2}{\partial z} \circ F^{-1}\right) \frac{\partial}{\partial p}
=:
g_1 \frac{\partial}{\partial s} 
+
g_2 \frac{\partial}{\partial p}.
\]
Therefore,
\[
\left.\nabla_{X}^{\mathbb{G}_2} Y\right\rvert_{(0,0)}
=
\left(\left. X(g_1) \frac{\partial}{\partial s} \right)\right\vert_{(0,0)}
+
\left(\left(\left. X(g_2)-\frac{4}{5} g_1^2\right) \frac{\partial}{\partial p} \right)\right\vert_{(0,0)}.
\]

\begin{prop}
The royal isometry $F_R:=(f_{R1},f_{R2}) : (\mathbb{D},ds^2_2) \rightarrow (\mathbb{G}_2, ds^2)$ defined by
\[
F_R(z) = (2z,z^2),
\quad z\in\mathbb{D}
\]
and the non-royal isometry $F_{N}:= (f_{N1}, f_{N2}): (\mathbb{D}, ds^2) \rightarrow (\mathbb{G}_2,ds^2)$ defined by
\[
F_N(z) = \left(
\sqrt{2} (\sqrt{z+1} - 1), 1-\sqrt{z+1}
\right),
\quad z\in\mathbb{D}.
\]
are not totally geodesic.
\end{prop}

\begin{proof}
We will first show that $F_{R}$ is not totally geodesic. From the previous discussion, 
\[
X = Y = (F_R)_*\left(\frac{\partial}{\partial z}\right) 
=
\left(\frac{\partial f_{R1}}{\partial z} \circ F^{-1}\right) \frac{\partial}{\partial s} 
+
\left(\frac{\partial f_{R2}}{\partial z} \circ F^{-1}\right) \frac{\partial}{\partial p}
=
2 \, \frac{\partial}{\partial s} 
+
s \, \frac{\partial}{\partial p}.
\]
This gives
\[
\left.\nabla_{X}^{\mathbb{G}_2} Y\right\rvert_{(0,0)}
=
\left(2 - \frac{4}{5} \cdot 4 \right) \left.\frac{\partial}{\partial p}\right\rvert_{(0,0)}
=
\frac{-6}{5} \left.\frac{\partial}{\partial p}\right\rvert_{(0,0)}.
\]
This implies that
\[
\left.\mathrm{II}(X,Y)\right\rvert_{(0,0)} \neq 0.
\]
Therefore, $F_R$ is not totally geodesic.

\medskip

Corresponding to the isometric embedding $F_N$, we have
\[
X = Y = (F_N)_*\left(\frac{\partial}{\partial z}\right) 
=
\left(\frac{\partial f_{N1}}{\partial z} \circ F^{-1}\right) \frac{\partial}{\partial s} 
+
\left(\frac{\partial f_{N2}}{\partial z} \circ F^{-1}\right) \frac{\partial}{\partial p}
=
\frac{1}{\sqrt{2} (1 - p)} \frac{\partial}{\partial s} 
-
\frac{1}{2 (1-p)} \frac{\partial}{\partial p}.
\]
This gives
\[
\left.\nabla_{X}^{\mathbb{G}_2} Y\right\rvert_{(0,0)}
=
-\frac{1}{2\sqrt{2}} \left.\frac{\partial}{\partial s}\right\rvert_{(0,0)}
+
\left(
\frac{1}{4} - \frac{2}{5}
\right) \left.\frac{\partial}{\partial p}\right\rvert_{(0,0)}
=
-\frac{1}{2\sqrt{2}} \left.\frac{\partial}{\partial s}\right\rvert_{(0,0)}
-
\frac{3}{20}
\left.\frac{\partial}{\partial p}\right\rvert_{(0,0)}.
\]
Thus, the second fundamental form does not vanish identically, and therefore $F_N$ is not totally geodesic.
\end{proof}


\section{Understanding the form of isometry} 

\subsection{k = 3 : Isometry with boundary values in $\mathbf{\mathcal{B}_3}$.}

We prove that the royal isometry is the only holomorphic isometric embedding upto automorphisms.

\begin{thm}\label{uni3}
Let $F=(S,P):(\mathbb{D},3 ds^2)\rightarrow(\mathbb{G}_2,ds^2)$ be a holomorphic isometric embedding. Then
\[
F(z) = (2z, z^2)
\]
after a possible reparametrization.
\end{thm}

\begin{proof}
Since, the cluster set of $\partial \mathbb{D}$ under $F$ lies in $\mathcal{B}_3$, we have $S^2 \equiv 4P$ on $\partial\mathbb{D}$. Then we must have by the identity principle that $S(z)^2= 4 P(z)$ for all $z\in\mathbb{D}$. Thus, $F(\mathbb{D}) \subset \mathcal{R}$. Since the automorphism group of $\mathbb{G}_2$ is transitive on $\mathcal{R}$, we can assume without loss of generality that $F(0) = (0,0)$. Therefore,
\[
K_{\mathbb{G}_2}((S(z),P(z)) = \frac{2}{\pi^2} \frac{1}{(1-\vert z\vert^2)^{6}}.
\]
This implies that
\[
\frac{16}{\pi^2} \left(\frac{1}{(2-\vert S(z)\vert^2 + 2\vert P(z)\vert^2)^3}\right) = \frac{2}{\pi^2} \frac{1}{(1-\vert z\vert^2)^{6}}.
\]
Since
\[
2-\vert S\vert^2 + 2\vert P\vert^2 = 2 - \vert S\vert^2 + \frac{\vert S\vert^4}{8} = 2\left(1-\frac{\vert S\vert^2}{4}\right)^2,
\]
we have
\[
\left(1 -\left(\frac{\vert S(z)\vert}{2}\right)^2 \right)^6 = (1-\vert z\vert ^2)^{6}.
\]
This implies that $\vert S(z)\vert = 2 \vert z\vert$ for all $z\in\mathbb{D}$. After a possible rotation, we have $S(z) = 2z$, and therefore, $P(z) = z^2$.
\end{proof}

\subsection{Isometries with lifts} 

Let $F:\mathbb{D}\rightarrow \mathbb{D}\times\mathbb{D}$ be a proper holomorphic map such that $F(0) = (0,0)$. We shall prove that the holomorphic map $\pi\circ F: (\mathbb{D}, k ds^2)\rightarrow (\mathbb{G}_2, ds^2)$ can never be an isometric embedding, for $k=1$ or $k=2$.

\begin{prop}\label{P'}
Let 
\[
F:= (S,P):(\mathbb{D},k ds^2)\rightarrow(\mathbb{G}_2,ds^2)
\]
be a holomorphic isometric embedding, where $k=1$ or $k=2$. Assume that $F(0) = (0,0)$. Then $P'(0) \neq 0$.
\end{prop}

\begin{proof}
Let $F:=(S, P)$. 
For $(s,p),(t,q)\in \mathbb{G}_2$, we know
\[
K_{\mathbb G_2}((s,p),(t,q))
=
\sum_{\alpha_1>\alpha_2\geq 0}
\frac{(\alpha_1+1)(\alpha_2+1)}{\pi^2}
f_\alpha(s,p)\overline{f_\alpha(t,q)}.
\]
We will first compute the first few $f_\alpha$'s.
\begin{eqnarray*}
f_{(1,0)}\circ\pi
&=&
\frac{z_1-z_2}{z_1-z_2}
=1\\
f_{(2,0)}\circ\pi
&=&
\frac{z_1^2-z_2^2}{z_1-z_2}
=
z_1+z_2
=s\\
f_{(2,1)}\circ\pi
&=&
\frac{z_1^2z_2-z_2^2z_1}{z_1-z_2}
=
\frac{z_1z_2(z_1-z_2)}{z_1-z_2}
=
z_1z_2
=p\\
f_{(3,0)}\circ\pi
&=&
\frac{z_1^3-z_2^3}{z_1-z_2}
=
z_1^2+z_1z_2+z_2^2
=
(z_1+z_2)^2-z_1z_2
=
s^2-p\\
f_{(3,1)}\circ\pi
&=&
\frac{z_1^3z_2-z_2^3z_1}{z_1-z_2}
=
z_1z_2
\frac{z_1^2-z_2^2}{z_1-z_2}
=
z_1z_2(z_1+z_2)
=
sp\\
f_{(3,2)}\circ\pi
&=&
\frac{z_1^3z_2^2-z_2^3z_1^2}{z_1-z_2}
=
z_1^2z_2^2
\frac{z_1-z_2}{z_1-z_2}
=
(z_1z_2)^2
=
p^2.
\end{eqnarray*}
Therefore,
\[
K_{G_2}((s,p),(t,q))
=
\frac{2}{\pi^2}
+
\frac{3}{\pi^2}s\overline{t}
+
\frac{6}{\pi^2}p\overline{q}
+
\frac{4}{\pi^2}(s^2-p)\overline{(t^2-q)}
+
\frac{8}{\pi^2}sp\,\overline{tq}
+
\frac{12}{\pi^2}p^2\overline{q}^2
+\cdots .
\]
The functional equation is
\begin{equation}\label{b1.1}
K_{\mathbb{G}_2}((S(z),P(z)),(S(w),P(w))) 
=
\frac{2}{\pi^2} \frac{1}{(1-z\ov{w})^{2k}}. 
\end{equation}
for $z,w\in\mathbb{D}$. The expansion of LHS in (\ref{b1.1}) gives
\begin{multline*}
K_{G_2}((S(z),P(z)),(S(w),P(w)))
=
\frac{2}{\pi^2}
+
\frac{3}{\pi^2}S(z)\overline{S(w)}
+
\frac{6}{\pi^2}P(z)\overline{P(w)}
\\+
\frac{4}{\pi^2}(S(z)^2-P(z))\overline{(S(w)^2-P(w))}
\\+
\frac{8}{\pi^2}S(z)P(z)\,\overline{S(w)P(w)}
+
\frac{12}{\pi^2}P(z)^2\overline{P(w)}^2
+\cdots .
\end{multline*}
The expansion of RHS in (\ref{b1.1}) for $k=1$ and $k=2$ gives
\[
\frac{2}{\pi^2}\frac{1}{(1-z\overline w)^2}
=
\frac{2}{\pi^2}
+
\frac{4}{\pi^2}z\overline w
+
\frac{6}{\pi^2}z^2\overline w^2
+\cdots 
\]
and
\[
\frac{2}{\pi^2}\frac{1}{(1-z\overline w)^4}
=
\frac{2}{\pi^2}
+\frac{8}{\pi^2}z\overline w
+\frac{20}{\pi^2}z^2\overline w^2
+\cdots ,
\]
respectively. Let, if possible, $P'(0) = 0$.
Write
\[
S(z)=c_1z+c_2z^2+\cdots 
\quad\text{and}\quad
P(z) = b_2 z^2 + b_3 z^3+\cdots.
\]
\begin{proof}[Case: k=1]
Comparing the coefficient of $\overline w$ in the identity (\ref{b1.1}), we get
\[
\frac{3}{\pi^2}S(z)\overline{c_1}
=
\frac{4}{\pi^2}z.
\]
Hence
\[
S(z)=Cz,
\quad
|C|^2=\frac43.
\]
Comparing the coefficient of $\overline w^2$ in (\ref{b1.1}), we obtain
\[
4(S(z)^2-P(z))\overline{(C^2-b_2)}
+
6P(z)\,\overline{b_2}
=
6z^2.
\] 
Substituting $S(z) = Cz$ gives
\[
\left(
-4\overline{(C^2-b_2)}
+
6\overline{b_2}
\right)P(z)
+
4C^2\overline{(C^2-b_2)}z^2
=
6z^2.
\]
We claim that the coefficient of $P(z)$ is nonzero. Indeed, if $-4\overline{(C^2-b_2)}
+
6\overline{b_2}
=
0$,
then 
\[
b_2=\frac25 C^2.
\]
So the coefficient of $z^2$ on the left-hand side would be
\[
4C^2\overline{\left(C^2-\frac25 C^2\right)}
=
4C^2\overline{\left(\frac35 C^2\right)}
=
\frac{12}{5}|C|^4.
\]
Using $|C|^2=\frac43$, this equals $\frac{64}{15}$, which is not equal to $6$. This is a contradiction. Thus, we have
\[
P(z)=Dz^2
\]
for some constant $D$. Substituting $S(z) = Cz$ and $P(z) = Dz^2$ in the diagonal functional equation gives
\[
\frac{16}{\pi^2}
\frac{
2-|C|^2|z|^2+2|D|^2|z|^4
}{
\left(
(2-|C|^2|z|^2+2|D|^2|z|^4)^2
-
|C^2-4D|^2|z|^4
\right)^2
}
=
\frac{2}{\pi^2}\frac{1}{(1-|z|^2)^2}.
\]
Put
\[
t=|z|^2,\qquad
A=|C|^2,\qquad
B=|D|^2,\qquad
E=|C^2-4D|^2.
\]
Then, we have
\[
\left((2-At+2Bt^2)^2-Et^2\right)^2
=
8(2-At+2Bt^2)(1-t)^2.
\]
The right-hand side is a
polynomial of degree at most $4$. However, on the left-hand side, the leading term is $16B^4t^8$. Therefore, $D = B=0$. But, $F(z) = (Cz, 0)$ is not an isometry for any value of $C$. Hence, a contradiction.
\end{proof}

\begin{proof}[Case: k=2] Comparing the coefficient of $\overline w$ in (\ref{b1.1}), we have
\[
3S(z)\overline{c_1}=8z.
\]
Therefore
\[
S(z)=Cz,
\quad
|C|^2=\frac{8}{3}.
\]
Comparing the coefficient of $\overline w^2$ in (\ref{b1.1}), and substituting $S(z)=Cz$, we get
\[
\bigl(10\overline{b_2}-4\overline{C}^2\bigr)P(z)
+
4C^2(\overline{C}^2-\overline{b_2})z^2
=
20z^2.
\]
If $10\overline{b_2}-4\overline{C}^2 = 0$, then $b_2=\frac{2}{5}C^2$.
So, the coefficient of $z^2$ on the left-hand side becomes $\frac{12}{5}|C|^4$.
Since $|C|^2=\frac{8}{3}$, this equals $\frac{256}{15}$, which is not equal to $20$. This is a contradiction. 
Therefore,
\[
P(z)=Dz^2
\]
for some constant $D$. Thus, isometry functional equation will give us the polynomial equation 
\[
\left((2-At+2Bt^2)^2-Et^2\right)^2
=
8(2-At+2Bt^2)(1-t)^4,
\]
where the notations are as before.
If $\vert D\vert^2 = B\neq 0$, then the left-hand side has degree $8$ in $t$, with leading term $16B^4t^8$,
whereas the right-hand side has degree at most $6$. Hence, $B=0$. This again forces $P\equiv 0$. Hence, a contradiction.
\end{proof}
Therefore, we must have $P'(0) \neq 0$. This finishes the proof.
\end{proof}

\begin{thm}\label{lift}
Let $F:=(f,g):\mathbb{D}\rightarrow\mathbb{D}\times\mathbb{D}$ be a proper holomorphic map with $F(0) = (0,0)$. Then, the holomorphic map $\pi\circ F : (\mathbb{D},k ds^2)\rightarrow(\mathbb{G}_2,ds^2)$ cannot be an isometric embedding either for $k=1$ or $k=2$. 
\end{thm}

\begin{proof}
Let $\pi\circ F(z) = (S(z), P(z))$. Here, $\pi\circ F(0) = (0,0)$. Since $P(z) = f(z) g(z)$ and the functions $f,g :\mathbb{D}\rightarrow\mathbb{D}$ vanish at the origin, we have $P'(0) = 0$. The result now follows by an application of Proposition \ref{P'}.
\end{proof}

\begin{cor}\label{nolift}
The non-royal isometry $F:=(S,P):(\mathbb{D},ds^2)\rightarrow (\mathbb{G}_2,ds^2)$ which is defined by
\[
F(z) := (S,P)= \left(\sqrt{2}(\sqrt{z+1} - 1), 1-\sqrt{z+1}\right),
\]
is a non-liftable isometric embedding.
\end{cor}

\begin{proof}
Note that $F(0) = (0,0)$. If there were a holomorphic map $G:=(g_1,g_2):\mathbb{D}\rightarrow \mathbb{D}\times\mathbb{D}$ such that $\pi\circ G = F$, then we must have $g_1(0) = g_2(0) = 0$. Since $F$ is proper, $G$ must also be a proper map. This will contradict Theorem \ref{lift}. Therefore, $F$ is non-liftable. 
\end{proof}

\subsection{No isometric embedding when k=2.}

We shall finally prove that the isometric constant $k$ cannot take the value $k=2$, that is, there does not exist any holomorphic isometric embedding $F: (\mathbb{D}, 2ds^2)\rightarrow (\mathbb{G}, ds^2)$. Therefore, the isometric constant can essentially take the values $k=1$ and $k=3$ only.

\medskip

We will first show that if such an isometry would exist, then its image must intersect the royal variety.

\begin{prop}\label{noroyal}
Let $F:(\mathbb D,2 ds^2)\rightarrow (\mathbb G_2,ds^2)$ be a proper holomorphic map such that
\[
F(\mathbb D)\cap \mathcal R=\varnothing.
\]
Then, $F$ is not an isometric embedding.
\end{prop}

\begin{proof}
Since $F(\mathbb{D})$ does not intersect the royal variety, $\mathbb D$ is simply connected, and
\[
    \pi:\mathbb D^2\setminus \Delta\rightarrow \mathbb G_2\setminus \mathcal R
\]
is a holomorphic covering, the function $F$ admits a
holomorphic lift
\[
H:\mathbb D\rightarrow \mathbb D^2\setminus \Delta
\]
satisfying $\pi\circ H=F$.
Moreover, since $F$ is proper and $\pi$ is proper, the lift $G$ is proper.
Let, if possible, $\pi\circ H$ be an isometry. Put $H: (f,g)$.
Since $F(\mathbb D)\cap \mathcal R=\varnothing$, we have
\[
f(z)\neq g(z),
\quad z\in\mathbb D.
\]
Without loss of generality, we can assume that $f(0) = 0$ and $g(0) = a\in \mathbb{D}\setminus\{0\}$.
Note that
\[
K_{\mathbb G_2}(F(z))
\approx \frac{1}{(1-|z|^2)^4}
\]
and $C(F)\subset \mathcal{B}_2$.
Here, $f$ and $g$ are
proper holomorphic self-maps of $\mathbb D$, and therefore finite Blaschke products. 
Since $f$ and $g$ are finite Blaschke products, they extend holomorphically past
$\partial\mathbb D$, and
\[
|f(\zeta)|=|g(\zeta)|=1
\qquad (\zeta\in\partial\mathbb D).
\]
We claim that there exists $\zeta_0\in\partial\mathbb D$ such that
\[
f(\zeta_0)=g(\zeta_0).
\]
Suppose not. Since $f(z)\neq g(z)$ for $z\in\mathbb D$, this would imply that $f-g$ has
no zeros on $\overline{\mathbb D}$. Hence, the index number of the curve $(f-g)(\partial\mathbb{D})$ around $0$ is
\[
n((f-g)(\partial\mathbb{D});0)
=
\frac{1}{2\pi i}\int_{(f-g)(\partial\mathbb{D})} \frac{1}{z} \,dz
=
0
\]
On the other hand, on $\partial\mathbb D$,
\begin{equation}\label{index}
f-g
=
g\left(\frac{f}{g}-1\right).
\end{equation}
Now, $\frac{f}{g}(\partial\mathbb{D})\subset \partial\mathbb{D}$ but it cannot be the entire $\partial\mathbb{D}$ as $f(z) \neq g(z)$ for all $z\in\ov{\mathbb{D}}$. Therefore, $\frac{f}{g}(\partial\mathbb{D})$ is a contractible closed curve contained in $\partial\mathbb{D}\setminus\{1\}$. Thus, $\gamma(z):=\frac{f}{g}(z) -1$ for $z\in\partial\mathbb{D}$, defines a contractible closed curve in $\mathbb{C}\setminus 0$. In terms of index number, we have
\[
n(\gamma;0) = 0.
\]
Thus, (\ref{index}) gives
\[
0 = n(g(\partial\mathbb{D});0) 
=
\frac{1}{2\pi i}\int_{g(\partial\mathbb{D})}\frac{1}{z} \,dz
=
\frac{1}{2\pi i}\int_{\partial\mathbb{D}}\frac{g'(z)}{g(z)} \,dz
=
\text{number of zeroes of }g \text{ in }\mathbb{D}
=
\deg(g).
\]
Therefore, $g$ must be a constant, giving a contradiction.
Hence, there exists $\zeta_0\in\partial\mathbb D$ such that
\[
f(\zeta_0)=g(\zeta_0).
\]
Let $B$ be a finite Blaschke product. Then
\[
B(z)=e^{i\theta}\prod_{j=1}^{N} b_{a_j}(z),
\qquad 
b_{a_j}(z)=\frac{z-a_j}{1-\overline{a_j}z},
\qquad a_j\in\mathbb D.
\]
For one Blaschke factor $b_{a_i}$, we have
\[
1-|b_{a_i}(z)|^2
=
1-\left|\frac{z-a_i}{1-\overline a_i z}\right|^2
=
\frac{(1-|a_i|^2)(1-|z|^2)}
{|1-\overline{a_i} z|^2}.
\]
Therefore, if we approach a boundary point $z_0\in \partial\mathbb{D}$ radially, we get
\[
1-|B(rz_0)|^2 = \sum_{k=1}^N (1-|b_{a_k}(rz_0)|^2)\prod_{j<k} |b_{a_j}(rz_0)|^2
\approx 1-r^2.
\]
Since $f$ and $g$ are finite Blaschke products, approaching $\zeta_0$ radially gives
\[
1-|f(r\zeta_0)|^2\approx 1-r^2,
\qquad
1-|g(r\zeta_0)|^2\approx 1-r^2.
\]
Moreover, since $f(\zeta_0)=g(\zeta_0)$, we have
\begin{eqnarray*}
|1-f(r\zeta_0)\overline{g(r\zeta_0)}|^2
&=&
(1-|f(r\zeta_0)|^2)(1-|g(r\zeta_0)|^2) + |f(r\zeta_0)-g(r\zeta_0)|^2.
\end{eqnarray*}
Since $(f-g)(\zeta_0)=0$, there exists a constant $C>0$ such that
\[
|f(r\zeta_0)-g(r\zeta_0)|^2
\leq C|r\zeta_0-\zeta_0|^2
=
C(1-r)^2\leq C(1-r^2)^2.
\]
As,
\[
|1-f(r\zeta_0)\overline{g(r\zeta_0)}|^2
\geq
(1-|f(r\zeta_0)|^2)(1-|g(r\zeta_0)|^2),
\]
we obtain
\[
|1-f(r\zeta_0)\overline{g(r\zeta_0)}|^2
\approx
(1-r^2)^2.
\]
Since
\[
K_{\mathbb G_2}(F(z))
=
\frac{(1-|f(z)|^2)(1-|g(z)|^2)+|1-f(z)\overline{g(z)}|^2}
{\pi^2 (1-|f(z)|^2)^2\,(1-|g(z)|^2)^2 \,|1-f(z)\overline{g(z)}|^4}.
\]
we must have
\[
K_{\mathbb G_2}(F(r\zeta_0))
\approx
\frac{1}{(1-r^2)^6}.
\]
This is not possible by Proposition \ref{isoc}.
Therefore, $F$ cannot be a holomorphic isometric embedding, which finishes the proof.
\end{proof}

We are now ready to prove our claim.

\begin{thm}\label{d=1}
Let
\[
F:=(S,P):(\mathbb D,2 ds^2)\rightarrow (\mathbb{G}_2,ds^2)
\]
be a holomorphic map. Then $F$ cannot be an isometric embedding.
\end{thm}

\begin{proof}
Let, if possible, there exists a holomorphic isometric embedding $F:=(S,P):(\mathbb D,2 ds^2)\rightarrow (\mathbb{G}_2,ds^2)$. It follows directly from Proposition \ref{noroyal} that $F(\mathbb{D})$ must intersect the royal variety $\mathcal{R}$. Since $Aut(\mathbb{G}_2)$ is transitive on $\mathbb{G}_2 \cap \mathcal{R}$, we can assume, without loss of generality, that $F(0) = (0,0)$. Moreover, we must have $P'(0) \neq 0$ from Proposition \ref{P'}.

\medskip

Moreover, Proposition \ref{cluster_set} says that $S$ and $P$ are rational functions. Let $N = \deg P$. Then, from the proof of Proposition \ref{cluster_set}, we can write
\[
P(z)=\tau z\frac{\widetilde D(z)}{D(z)},
\qquad |\tau|=1,
\]
where
\[
D(z)=\prod_{j=1}^{N-1}(1-\overline{\alpha_j}z),
\qquad
\widetilde D(z)=z^{N-1}\overline{D(1/\overline z)},
\qquad
\alpha_j\in\mathbb{D}\setminus\{0\}.
\]
This is because $P(0)= 0$ and $P'(0)\neq 0$.
Let $d$ denote the leading coefficient of $D$, which is $(-1)^{N-1}\prod_{j=1}^{N-1}\ov{\alpha_j} \neq 0$. Then $\widetilde{D}$ has the constant term equal to $\ov{d}$. Thus, $P'(0) = \tau \ov{d}$. Furthermore,
\[
S(z) = \frac{Q(z)}{D(z)}
\]
where $Q$ is a polynomial of degree at most $N$. For $\zeta\in\mathbb{D}$, the function $Q$ satisfies
\begin{equation}\label{Q}
Q(\zeta) = \tau \zeta^N \overline{Q(\zeta)},
\end{equation}
Write
\[
Q(z) = \sum_{k=0}^{N} q_k z^k.
\]
Comparing the Fourier coefficients on both sides of (\ref{Q}) gives $q_0 = \tau \ov{q_N}$. Since $q_0 = Q(0) = S(0)= 0$, we have $\deg Q \leq N-1$.

\medskip

Using the expansion of the Bergman kernel of $\mathbb G_2$, we get
\begin{multline*}
K_{G_2}((S(z),P(z)),(S(w),P(w)))
=
\frac{2}{\pi^2}
+
\frac{3}{\pi^2}S(z)\overline{S(w)}
+
\frac{6}{\pi^2}P(z)\overline{P(w)}
\\+
\frac{4}{\pi^2}(S(z)^2-P(z))\overline{(S(w)^2-P(w))}
\\+
\frac{8}{\pi^2}S(z)P(z)\,\overline{S(w)P(w)}
+
\frac{12}{\pi^2}P(z)^2\overline{P(w)}^2
+\cdots .
\end{multline*}
On the other hand,
\[
\frac{2}{\pi^2}\frac{1}{(1-z\overline w)^4}
=
\frac{2}{\pi^2}
+\frac{8}{\pi^2}z\overline w
+\frac{20}{\pi^2}z^2\overline w^2
+\cdots .
\]
Write
\[
S(z)=c_1z+c_2z^2+\cdots
\quad\text{and}\quad
P(z)=b_1z + b_2z^2+b_3z^3+\cdots .
\]
Comparing the coefficient of $\overline{w}$ on both the sides of
\[
K_{\mathbb{G}_2}((S(z),P(z)),(S(w),P(w))) 
=
\frac{2}{\pi^2} \frac{1}{(1-z\ov{w})^{4}},
\]
we obtain
\begin{equation}\label{ccm}
3 \overline{c_1}S(z) + 10 \overline{b_1} P(z) - 4 \ov{b_1} S(z)^2  = 8 z. 
\end{equation}
Substituting
\[
S=\frac{Q}{D},
\qquad
P=\tau z\frac{\widetilde D}{D}
\]
in the equation (\ref{ccm}), and multiplying by $D^2$ on both sides, we get
\begin{equation}\label{ccm2}
3\overline{c_1}\,QD
+
10\overline{b_1}\,\tau z\widetilde D D
-
4\overline{b_1}\,Q^2
=
8z\,D^2. 
\end{equation}
We will compare the coefficient of $z^{2N-1}$ on both sides of (\ref{ccm2}).
Note that
\[
\deg(QD)\leq 2N-2,
\qquad
\deg(Q^2)\leq 2N-2.
\]
The coefficient of $z^{N-1}$ in $\widetilde{D}$ is $1$, and the coefficient of $z^{N-1}$ in $D$ is $d$. Therefore,
\[
10 \ov{b_1} \tau d = 8 d^2.
\]
Since $b_1 = P'(0) =  \tau \ov{d}$, we have
\[
10 d^2 = 8 d^2
\]
which is not possible unless $d=0$. Hence, a contradiction. Therefore, $F$ cannot be an isometric embedding.
\end{proof}

\section{Uniqueness of non-royal embedding}

We know from Theorem \ref{lift} that there does not exist any holomorphic isometric embedding $F_N:(\mathbb{D}, ds^2) \rightarrow (\mathbb{G}_2, ds^2)$ such that $F_N$ is liftable and $F(\mathbb{D})$ intersects the royal variety. We have an example of an isometric embedding $F:(\mathbb{D}, ds^2) \rightarrow (\mathbb{G}_2, ds^2)$ in Example \ref{nonroyal}, defined by
\[
F_N(z) := (S,P)= \left(\sqrt{2}(\sqrt{z+1} - 1), 1-\sqrt{z+1}\right).
\]

The image $F_N(\mathbb{D})$ intersects the royal variety at exactly one point $(0,0)$. The isometry $F_N$ is non-liftable, and $C(F_N) \subset \mathcal{B}_1 \cup \mathcal B_2$ with $C(F_N)\cap \mathcal{B}_2 =\{ (-2,1)\}$. In this section, we shall prove that $F_N$ is the unique isometric embedding upto a reparametrization.

\medskip

Let $F:= (S(z),P(z)):(\mathbb{D}, ds^2)\rightarrow (\mathbb{G}_2, ds^2)$ be a holomorphic isometric embedding with $F(0) = x_0 = (s_0,p_0)$.

\begin{lem}
There exists a rational function $R(s,p)$ such that $R(F(z))=z$ for every $z\in\D$.
\end{lem}

\begin{proof}
We know that the Bergman kernel $K_{\mathbb{G}_2}$ of the symmetrized bidisk is rational. For $(s,p)\in\mathbb{G}_2$ and $z\in\mathbb{D}$, define the rational function
\[
\mathcal K((s,p),w)=\frac{K_{\G}((s,p),F(w))K_{\G}((s_0,p_0),(s_0,p_0))}{K_{\G}((s,p),(s_0,p_0))K_{\G}((s_0,p_0),F(w))},
\]
and
\[
R((s,p))=\frac12\left.\frac{\partial\mathcal K((s,p),w)}{\partial\overline w}\right|_{w=0}.
\]
It is immediate that $R$ is a rational function in $(s,p)\in\mathbb{C}^2$. The functional equation of the isometry gives
\begin{equation}\label{feq}
\frac{K_{\G}(F(z),F(w))K_{\G}((s_0,p_0),(s_0,p_0))}
{K_{\G}(F(z),(s_0,p_0))K_{\G}((s_0,p_0),F(w))}
=\frac1{(1-z\overline w)^2}.
\end{equation}
Therefore,
\[
R(F(z))=\frac12 \left.\frac{\partial}{\partial{\ov w}}\frac{1}{(1-z\ov w)^{2}}\right\rvert_{w=0}
=z
\]
for all $z\in\mathbb{D}$.
\end{proof}

\begin{prop}
There is a nonzero irreducible polynomial $G(s,p)\in\mathbb C[s,p]$ such that
\begin{equation}
G(S(z),P(z))=0,
 \qquad z\in\D. 
\end{equation}
Equivalently, $F(\mathbb{D})$ is contained in an irreducible projective
algebraic curve $\Gamma\subset\mathbb P^2$.
\end{prop}

\begin{proof}
Assume that no such nonzero polynomial exists. This means that
$F(\D)$ is \emph{Zariski dense} in $\mathbb C^2$, that is, a polynomial vanishing on $F(\mathbb{D})$ must be the zero polynomial. For $z,w\in\mathbb{D}$, we have
\[
\frac{K_{\G}(F(z),F(w))K_{\G}((s_0,p_0),(s_0,p_0))}
{K_{\G}(F(z),(s_0,p_0))K_{\G}((s_0,p_0),F(w))}
=
\frac1{(1-R(F(z))\overline{R(F(w))})^2}.
\]
After clearing the denominators, this is a polynomial equation in the four
independent variables $S(z),P(z),\overline{S(w)},\overline{P(w)}$. 
The product of two Zariski dense subsets of $\mathbb C^2$ is Zariski dense in
$\mathbb C^4$, this can be seen by restricting the polynomial to each of the Zariski dense sets, successively. Since $\overline{F(\mathbb{D})}$ is also Zariski dense in $\mathbb{C}^2$, the polynomial equation is identically zero. Therefore,
\[
 \frac{K_{\G}((s,p),(t,q))K_{\G}((s_0,p_0),(s_0,p_0))}
 {K_{\G}((s,p),(s_0,p_0))K_{\G}((s_0,p_0),(t,q))}
 =\frac1{(1-R((s,p))\overline{R((t,q))})^2}
\]
for arbitrary $(s,p),(t,q)\in\mathbb{G}_2$ near $(s_0,p_0)$.
Thus
\[
 ds^2_{\G}=R^*ds^2_\D\qquad\text{near }(s_0,p_0).
\]
The left-hand side is a positive-definite Hermitian form of complex rank two,
while the right-hand side is the pullback of the Bergman metric on a one-dimensional complex
manifold and has rank at most one. Contradiction.
Hence, there exists a nonzero irreducible polynomial $G(s,p)\in \mathbb{C}[s,p]$ such that
\[
G(S(z),P(z))=0,
 \qquad z\in\D. 
\]
Homogenize the polynomial $G$ to obtain an irreducible projective algebraic curve $\Gamma\subset\mathbb{P}^2$. 
\end{proof}

Since $\Gamma$ is a plane projective algebraic curve (projective algebraic hypersurface in $\mathbb{P}^2$), there exists a (unique up to biholomorphisms) normalization $\sigma: \widetilde{\Gamma} \rightarrow \Gamma$ of $\Gamma$, where $\widetilde{\Gamma}$ is a compact Riemann surface, and $\sigma$ is holomorphic; see \cite[\S 6, Chapter II]{Griffiths}. Since $F:\mathbb{D}\rightarrow \Gamma$ is a non-constant holomorphic map, there exists a holomorphic lift $\widetilde F : \mathbb{D}\rightarrow \widetilde\Gamma$ satisfying $\sigma \circ \widetilde{F} = F$.
Since $R(s,p)$ is rational for $(s,p)\in\mathbb{C}^2$, we have the meromorphic function $\widetilde{R} := R\circ \sigma : \widetilde{\Gamma} \rightarrow \mathbb{P}^1$ on $\widetilde{\Gamma}$. For $z\in\mathbb{D}$, we have
\[
\widetilde{R}(\widetilde{F}(z))
=
(R\circ \sigma)(\widetilde{F}(z))
=
R(F(z)) = z.
\]

\medskip

Now, for $X = (t, s, p)$ and $Y = (t', s',p')$ in $\mathbb{C}^3$, define
\[
\begin{aligned}
 H(X,Y)&=2t\overline{t'}-s\overline{s'}+2p\overline{p'},
 & Q(X)&=s^2-4tp,\\
 D(X,Y)&=H(X,Y)^2-Q(X)\overline{Q(Y)}.
\end{aligned}
\]
The Hermitian form $H$ has matrix $J=\operatorname{diag}(2,-1,2)$ and is
\emph{nondegenerate}: if $H(X,Y)=0$ for every $Y$, then $X=0$.
Write $\widehat x=(1,s,p)$ for $x=(s,p)\in\mathbb{C}^2$.
Therefore, the Bergman kernel of $\mathbb{G}_2$ can be written as
\begin{equation}\label{bergman}
 K_{\G}(x,y)=\frac{16}{\pi^2}\,
 \frac{H(\widehat x,\widehat y)}{D(\widehat x,\widehat y)^2},
 \qquad
 x,y\in\mathbb{G}_2.
\end{equation}
From \eqref{bergman}, note that we can talk about $K_{\mathbb{G}_2}(x,y)$ for all $x,y\in\mathbb{P}^2$ by substituting $\widehat{x}, \widehat{y}$ with $x,y$ in the expression. For $w\in\D$, put $X_w=(1,S(w),P(w))$ and define
\[
 N_w=\{[X]\in\mathbb P^2:H(X,X_w)=0\}.
\]
Here, $N_w$ is a projective line because $H(\cdot, X_w)$
is a nonzero homogeneous linear form. Note that $H(X_w,X_w)>0$, so $N_w$
never contains the entire curve $\Gamma$. 

\begin{prop}\label{AI_help}
Let $\Gamma\subset\mathbb P^2$ be the irreducible projective
algebraic curve containing $F(\mathbb D)$.
Then $\Gamma$ is a projective line.
\end{prop}

\begin{proof}
The functional equation for the isometry $F$ gives 
\[
K_{\mathbb G_2}(F(z),F(w))\,K_{\mathbb G_2}(x_0,x_0)\,
(1-z\overline w)^2
=
K_{\mathbb G_2}(F(z),x_0)\,
K_{\mathbb G_2}(x_0,F(w))
\]
for $z,w\in\mathbb{D}$. 
This can be written as
\begin{equation}\label{fe}
K_{\mathbb G_2}(\sigma \circ \widetilde{F}(z),F(w))\,K_{\mathbb G_2}(x_0,x_0)\,
(1-z\overline w)^2
=
K_{\mathbb G_2}(\sigma \circ \widetilde{F}(z),x_0)\,
K_{\mathbb G_2}(x_0,F(w)).
\end{equation}
Since $\widetilde{R}(\widetilde{F}(z)) = z$ for all $z\in \mathbb{D}$ and $\widetilde{F}(\mathbb{D})$ is Zariski dense in $\widetilde{\Gamma}$, the identity \eqref{fe} extends meromorphically to $\widetilde{\Gamma}$. Therefore, for each $w\in\mathbb{D}$, 
\begin{equation}\label{fe_extension}
K_{\mathbb G_2}(\sigma(x),F(w))\, K_{\mathbb G_2}(x_0,x_0)\,
(1-\widetilde{R}(x)\,\overline w)^2
=
K_{\mathbb G_2}(\sigma(x),x_0)\,
K_{\mathbb G_2}(x_0,F(w)),
\end{equation}
where $x\in\widetilde{\Gamma}$.
Suppose that $\Gamma$ is not a line. Then, the vectors $X_w=(1,S(w),P(w))$ span $\mathbb C^3$, otherwise $F(\D)$ would lie in a projective line. 

\medskip

We will show that $Q$ vanishes identically on $\Gamma$. Suppose not. Let $E\subset\widetilde\Gamma$ consist of the points lying over the singular points of $\Gamma$ and the points at infinity, the zeros of $Q\circ\sigma$, the zeros and poles of $K_{\mathbb G_2}(\sigma(x),x_0)$, and the poles of $\widetilde{R}$. It is a standard fact that $\Gamma$ has only finitely many singularities, see \cite{Fischer}. Since the zeros or poles of
nonzero meromorphic functions on the compact Riemann surface $\widetilde\Gamma$ must be finite, we see that $E$ is finite.

\medskip

For each $\xi\in E$, the function $w\mapsto H(\sigma(\xi),X_w)$ is not identically zero for $w\in\mathbb{D}$. Indeed, otherwise the nondegeneracy of $H$ and the fact that the vectors $X_w$ span $\mathbb C^3$ would force $\sigma(\xi)=0$, which is impossible for a projective point. 
Moreover, $w\mapsto Q(X_w)$ is not identically zero for $w\in\mathbb{D}$, because $Q$ is not identically zero on $\Gamma$. We may therefore choose $w\in\mathbb D$ such that
\[
Q(X_w)\ne0
\quad\text{and}\quad
H(\sigma(\xi),X_w)\ne0
\quad\text{for every }\xi\in E.
\]
By the B\'ezout's theorem (see \cite[Section 2.7]{Fischer}), the projective line $N_w$ intersects
$\Gamma$. By our choice of $w$, none of these intersection points is in $\sigma(E)$.
Choose $y\in N_w\cap\Gamma$. Then $y$ is a smooth affine point, and
\[
Q(\widehat y)\ne0,\qquad
K_{\mathbb G_2}(y,x_0)\ne0,\qquad
R(y)\ne\infty.
\]
Since $H(\widehat y,X_w)=0$, we have
\[
D(\widehat y,X_w)
=-Q(\widehat y)\,\overline{Q(X_w)}\ne0.
\]
Consequently, $K_{\mathbb G_2}(y,F(w))=0$. Substituting $\sigma(x) = y$ into \eqref{fe_extension}, we obtain
\[
0=K_{\mathbb G_2}(y,x_0)\, K_{\mathbb G_2}(x_0,F(w)).
\]
If $K_{\mathbb G_2}(x_0,F(w)) = 0$, then $K_{\mathbb{G}_2}(F(z),F(w))= 0$ for all $z\in\mathbb{D}$, which is not possible. Therefore, we arrive at a contradiction. Hence, $Q$ vanishes identically on $\Gamma$. 

\medskip

In particular, we have proved that $S(z)^2 - 4 P(z) =0$ for all $z\in\mathbb{D}$.
Thus, $F(\mathbb{D})\subset\mathcal{R}$. However, in such a case, $F$ is an isometry only when $k=3$. Therefore, our assumption that $\Gamma$ is not a projective line was wrong. This proves the proposition.
\end{proof}

\begin{thm}
Let $F:=(S(z),P(z)):(\mathbb{D}, ds^2)\rightarrow (\mathbb{G}_2, ds^2)$ be a holomorphic isometric embedding. Then
\[
F(z) =\left(\sqrt{2}(\sqrt{z+1} - 1), 1-\sqrt{z+1}\right).
\]
after a possible reparametrization.
\end{thm}

\begin{proof}
We have proved that the irreducible projective algebraic curve $\Gamma$ containing $F(\mathbb{D})$ must be a projective line. Therefore, there exists $a,b\in\mathbb{C}$, not both zero, such that
\begin{equation}\label{final1}
   S(z) + b \,P(z) = a 
\end{equation}
for all $z\in \mathbb{D}$. Let $W = \operatorname{span}_{\mathbb{C}}\{ X_w : w\in \mathbb{D} \}$. Here, $W$ is a two-dimensional subspace of $\mathbb{C}^3$. We will first prove that $H\vert_{W}$ has rank one. 

\medskip

Let, if possible, $H\vert_W$ be nondegenerate, that is, it has rank two. By the B\'ezout's theorem, the projective lines $N_0$ and $\Gamma$ intersect at a single point with multiplicity $1$. Thus, there exists $\xi\in\Gamma$ such that $H(\xi, X_0) = 0$ with order $1$. 
Since $R(F(z)) = z$ for all $z\in\mathbb{D}$, the functional equation gives
\[
 \frac{H(X_z, X_w)}{H(X_z,X_0)}
 =
 \frac{K_{\G}(x_0,F(w))}{K_{\G}(x_0,x_0)}\left(\frac{D(X_z,X_w)}{D(X_z,X_0)(1-R(X_z) \,\ov w )}\right)^2,
 \qquad
 z,w\in\mathbb{D}.
\]
Moving to the normalization $\widetilde{\Gamma}$, we have the identity on $\widetilde F(\mathbb{D})$
\[
 \left(\frac{H(\sigma(X), X_w)}{H(\sigma(X),X_0)}\right) 
 =
 c_w\left(\frac{D(\sigma(X), X_w)}{H(\sigma(X), X_0)\, (1-\widetilde{R} \,\ov w )}\right)^2,
\qquad
c_w>0.
\]
Both sides are meromorphic functions on $\widetilde\Gamma$. 
The left hand side is equal to a nonzero constant times the square of a meromorphic
function on $\widetilde{\Gamma}$. Therefore,
\[
\operatorname{ord}_{P} H(\cdot, X_w) - \operatorname{ord}_{P} H(\cdot, X_0) \in 2\mathbb{Z}
\]
for every $w\in \mathbb{D}$ and $P\in\Gamma$. Since $\operatorname{ord}_{\xi} H(\cdot, X_0) = 1$, we have $\operatorname{ord}_{\xi} H(\cdot, X_w) = 1$ for every $w\in \mathbb{D}$. That is, $H(\xi,X_w) = 0$ for all $w\in\mathbb{D}$. Thus, $H\vert_W$ is degenerate. Since $H(X_w,X_w)>0$ for every $w\in\mathbb{D}$, the rank of $H\vert_W$ cannot be zero. Hence, we have proved that $H\vert_W$ has rank one.

\medskip

For $z\in\mathbb{D}$, we have 
\[
X_w = (1, S(w), P(w)) = (1, a - b\,P(w), P(w)) 
=
(1,a,0) + P(w)\,(0, -b, 1).
\]
With respect to the basis $\{(1,a,0),(0,-b,1)\}$ of $W$, the matrix of $H\vert_W$ is given by
\[
\begin{pmatrix}
2-|a|^2 & a\overline{b}\\
\overline{a}b & 2-|b|^2
\end{pmatrix}.
\]
Since the rank of $H\vert_W$ is one, its determinant vanishes. This gives
\begin{equation}\label{final2}
    \vert a\vert^2 + \vert b\vert^2 = 2.
\end{equation}
Now for each $z\in \mathbb{D}$, we have
\begin{eqnarray*}
H(X_z,X_z)
&=&
2-|a-bP(z)|^2+2|P(z)|^2
\\&=&
\vert b\vert^2 + \vert a \vert^2 \vert P(z)\vert^2 + a \,\overline{b P(z)} + \overline{a}\, bP(z) 
\\&=&
\vert b +a P(z)\vert^2.
\end{eqnarray*}
Put $q(z)=\overline{b}+\overline{a}P(z)$. Since $H(X_z,X_z)>0$, the holomorphic function
$q$ is nowhere vanishing on $\mathbb D$. Now define
\[
\varphi(z)=\frac{S(z)^2-4P(z)}{q(z)^2}.
\]
Since $D(X_z,X_z)=H(X_z,X_z)^2-|S(z)^2-4P(z)|^2>0$,
we have $|\varphi(z)|<1$ for every $z\in\mathbb D$. Thus, we have the holomorphic function $\varphi:\mathbb D\to\mathbb D$.
Moreover, the Bergman kernel formula of $\mathbb{G}_2$ gives
\begin{eqnarray*}
K_{\mathbb G_2}(F(z),F(z))
&=&
\frac{16}{\pi^2}\frac{\vert q(z)\vert^2}{(\vert q(z)\vert^4- \vert S(z)^2 - 4P(z)\vert^2)^2}
\\&=&
\frac{16}{\pi^2}\frac{1}{|q(z)|^6(1-|\varphi(z)|^2)^2}.
\end{eqnarray*}
Since $q$ is holomorphic and nowhere vanishing, applying $\partial\bar\partial\log$ on both sides yields
\[
ds_{\mathbb D}^2= F^*ds_{\mathbb G_2}^2
=\varphi^*ds_{\mathbb D}^2.
\]
The first equality holds because $F$ is an isometry. Thus, the holomorphic function $\varphi:(\mathbb{D},ds^2)\rightarrow(\mathbb{D},ds^2)$ is an isometry. This forces $\varphi$ to be an automorphism of $\mathbb{D}$. Hence, there exists $z_0\in\mathbb{D}$ such that $\varphi(z_0) = 0$. 

\medskip

We have proved that $S(z_0)^2-4P(z_0) = 0$. That is, $F(\mathbb{D})$ intersects the royal variety. Therefore, we can assume, after a reparametrization, that $F(0) = (0,0)$. Equations \eqref{final1} and \eqref{final2} gives
\begin{equation}
S(z) = -c \sqrt{2} P(z),
\qquad 
\vert c\vert =1.
\end{equation}
Rotating both roots of the symmetrization simultaneously by $e^{i\theta}$ changes $(s,p)$ to $(e^{i\theta} s, e^{2i\theta}p)$ and changes $c$ to $e^{-i\theta}c$. Therefore, we can assume, after a reparametrization, that $c=1$.
Finally, the functional equation gives
\begin{eqnarray*}
\frac{2}{\pi^2}\left(
\frac{1}{(1-\vert z\vert^2 )^2}
\right)
&=&
\frac{16}{\pi^2}\left(
\frac{2 - \vert S(z)\vert^2 + 2\vert P(z)\vert^2}{((2-\vert S(z)\vert^2 + 2\vert P(z)\vert^2)^2-\vert S(z)^2 - 4P(z)\vert^2)^2}
\right)
\\&=&
\frac{16}{\pi^2}\left(
\frac{2}{(4-\vert 2 P(z)^2 - 4P(z)\vert^2)^2}
\right)
\end{eqnarray*}
which can we rewritten as $4 - \vert 2 P(z)^2 - 4 P(z)\vert^2 = 4 - 4\vert z\vert^2$. So, $\vert P(z)^2 - 2 P(z)\vert = \vert z\vert$ for all $z\in\mathbb{D}$. Thus, there exists a unimodular constant $\eta$ such that
\[
P(z)^2 - 2 P(z) = \eta z.
\]
Solving the equation gives $P(z) = 1 -  \sqrt{1-\eta z}$ because $P(0) = 0$. The reparametrization $F(-\overline{\eta} \,z)$ finally gives
\[
P(z) = 1 - \sqrt{1+z} 
\qquad\text{and}\qquad
S(z) = \sqrt{2} ( \sqrt{1 + z} -1).
\]
This proves the uniqueness.
\end{proof}

\medskip

\textbf{Acknowledgment.} The author is grateful to Prof. Włodzimierz Zwonek for suggesting the problem and for the initial discussions that motivated this work.
This research was commenced when the author was a postdoctoral researcher at Jagiellonian University, Kraków, Poland, supported by SHENG-3 K/NCN/000286, research project UMO2023/Q/ST1/00048. The author also acknowledges the subsequent support from the ANRF-National Postdoctoral Fellowship, project no. PDF/2025/001158, at IISER Pune, India.

\medskip

\textbf{AI Declaration.} 
ChatGPT was used to suggest the application of Bézout's theorem to the projective lines $N_w$ as an approach to proving the results
of Section 6. ChatGPT was also used for language polishing. All mathematical statements and all other proofs were developed independently by the author. The author takes full responsibility for the content of this article.


\end{document}